\documentclass[10pt]{amsart}
\usepackage{amsfonts}
\usepackage{amssymb}
\usepackage{amsthm}
\usepackage[dvips]{graphicx}
\theoremstyle{plain}
\newtheorem{theorem}{Theorem}
\newtheorem{proposition}[theorem]{Proposition}
\newtheorem{lemma}[theorem]{Lemma}
\newtheorem{definition}[theorem]{Definition}
\newtheorem{corollary}[theorem]{Corollary}
\newtheorem{remark}[theorem]{Remark}

\input xy
\xyoption{all}

\newcommand{\dxy}{\frac{dx}{y}}
\newcommand{\dxys}{(\frac{dx}{y})^2}
\renewcommand{\epsilon}{\varepsilon}

\newcommand{\cbG}{\check{\bar{\Gamma}}}
\newcommand{\surj}{\twoheadrightarrow}

\renewcommand{\SS}{\mathcal{S}}
\newcommand{\LL}{\mathcal{L}}
\newcommand{\HH}{\mathcal{H}}
\newcommand{\WW}{\mathcal{W}}
\newcommand{\BB}{\mathcal{B}}
\newcommand{\MM}{\mathcal{M}}
\newcommand{\Tau}{T}

\begin{document}
\title{Monodromy of The Hitchin Map over Hyperelliptic Curves}
\author{D. Jeremy Copeland}
\maketitle

\section{introduction}
The purpose of this paper is to study the monodromy of the Hitchin fibration for rank 2 vector bundles over hyperelliptic curves.  We reduce the problem to studying a surface braid group generalization of the classical Burau representation \cite{Birman}, and give a combinatorial method for computing this representation.

Let $M$ be a Riemann surface of genus $g>1$ with canonical bundle $K$.  In  \cite{Hitchin}, Hitchin considered the set of stable Higgs pairs, $(V,\Phi)$, where $V$ is a rank 2 vector bundle of fixed odd degree and fixed determinant $\bigwedge^2V=\xi$, and $\Phi\in \text{End}(V)\otimes K$.  $(V,\Phi)$ is called stable if any $\Phi$-invariant line sub-bundle $L$ of $V$ has the property that $\text{deg}(L)<\frac12\text{deg}(V)$.  If  $\mathcal{G}$ is the group of automorphisms of $V$ with determinant one, then we study the space $\MM$ which is the quotient of the set of all stable Higgs pairs $(V,\Phi)$ modulo $\mathcal{G}$.  This is a manifold by results of Narasimhan and Ramanan \cite{Narasimhan}.  

Hitchin considered the map $\MM\to H^0(M,K^2), (V,\Phi)\mapsto \det(\Phi)$ and showed that the generic fiber is  an abelian variety.  Specifically, if we let $H^s\subset H^0(M,K^2)$ be the space of quadratic differentials with simple zeroes, then the restriction $\MM|_{H^s}$ is a fiber bundle of abelian varieties.  Monodromy of this bundle is the focus of this paper.

We consider the fibers explicitly.  Let $\omega\in H^s$.  There is a spectral cover $S=S_\omega\to M$ which is a ramified double cover of $M$ constructed as follows.  In the total space, $TK$, of the canonical bundle, $K$, take the curve $S=\sqrt{\omega}$ (Hitchin also shows that this is nonsingular).  This is the set of points which map to $\omega$ under the squaring map taking $TK$ to the total space of $K^2$.  Notice that $c:S\to M$ comes equipped with an involution, $\tau$, which fixes $M$ and swaps the sheets of the cover.  The Prym variety of $S$ is the subgroup of the Jacobian of $S$ on which $\tau$ acts by $-1$.  In our case, this is just $-1\in\mathbb C$ acting on $S\subset TK\to M$.  It is known from \cite{Hitchin} that if (and only if) the divisor associated to $\omega\in H^0(M,K^2)$ has simple zeroes, then the fiber of $\MM$ over $\omega$ is isomorphic to the Prym variety of $(S,M)$.
For a more complete description and reference set for Hitchin's system, see Section 6 of \cite{Hitchin} or Section 8 of \cite{Hitchin2}.  We also describe this in more detail in the proof of Theorem \ref{thm.MtoJ}.

One may also construct a bundle, $\SS\to H^s$ whose fibers are the spectral covers $S_\omega\subset TK_M$ and the associated bundle of Jacobians (Section \ref{sec.bundleS}), and ask whether this contains $\MM$ as $\text{Jac}(\SS)^{\tau=-1}$.  In the proof of Theorem \ref{thm.MtoJ}, we don't prove this, but we show that it is true up to translation, akin to the statement that $\text{Jac}(S)\cong\text{Pic}^k(S)$ for $k\neq 0$, due to a lack of canonical basepoint.  Because of this, these bundles agree homologically, since homology of an abelian variety is translation invariant.  Thus to study the monodromy of $\MM$, we may instead study the monodromy of the Jacobian varieties, or indeed of the $\tau=-1$ part of the homology of $\SS$.

In Section \ref{sec.relating}, we show that the $\tau=-1$ part of the homology of $\SS$ is the restriction of a generalization of the classical Burau representation to (a subgroup of) the surface braid group.  We then study this representation.

Throughout this paper, we consider only hyperelliptic curves, $M$.  This makes many of the arguments more concrete, but Theorem \ref{thm.iso} is the main point where this is used.  In Section \ref{sec.generalize}, we show how Theorem \ref{thm.iso} is in fact the only obstruction to generalizing.

The main result of this paper may be stated as follows:

\begin{theorem}\label{thm.teaser}
To each hyperelliptic $M$ of genus greater than 2, one may associate a graph $\check\Gamma$ with edge set $E$ and skew bilinear pairing $(e\cdot e')$ on $e,e'\in E$ such that
\begin{itemize}
\item The monodromy representation of $\pi_1(H^s)$ is generated by elements $\sigma_e$ labelled by the edges $e\in E$ and the element $\tau$.
\item The monodromy representation is a quotient of $\mathbb ZE$.
\item The action of $\pi_1(H^s)$ on $\mathbb ZE$ is given by:
\begin{align}
\tau e&=-e\notag\\
\sigma_{e_\alpha}e&=e-(e\cdot e_\alpha)e_\alpha.\notag
\end{align}
\end{itemize}
\end{theorem}
These results are, respectively,  Corollary \ref{cor.edgesgenerate} and Corollary \ref{cor.1tau}, Theorem \ref{thm.prymef}, and Theorem \ref{thm.repwithint}.  The intersection numbers are found in Proposition \ref{prop.intnos}, and Proposition \ref{prop.ZF} gives the kernel of the quotient.

\section{The fundamental group of $H^s$}\label{sec.geomresults}

In this paper, we study the monodromy action of $\pi_1(H^s)$ on this bundle of abelian varieties, $\MM|_{H^s}$.  We begin by discussing the nature of $\pi_1(H^s)$, and reduce to a study of braid groups.  In Proposition \ref{prop.withdiagram} we show that, up to $\tau$, the monodromy representation factors through a surface braid group, so that we need only understand the image of $\pi_1(H^s)$ in this group.  This image is the result of Theorem \ref{thm.iso}, for which we begin to lay the groundwork here.

Notice that if $q_1$ and $q_2$ are elements of $H^s$ with the same zeroes, then $q_1/q_2$ is constant, thus $H^s/{\mathbb{C}^*}$ is the space of effective divisors with simple zeroes linearly equivalent to $K^2$.  Call this space $PH^s$.  An element of $PH^s$ is determined exactly by its divisor.

\begin{proposition}
The kernel of the map induced from projection $\pi_1(H^s)\to\pi_1(PH^s)$ acts by $\{1,\tau\}$ in the monodromy.
\end{proposition}

\begin{proof}
As a quotient by a free $\mathbb C^*$ action, the map $H^s\to PH^s$ is a fibration with fiber $\mathbb{C}^*$, so
\[
\pi_1(\mathbb{C}^*)\to \pi_1(H^s)\twoheadrightarrow \pi_1(PH^s)
\]
is exact.  Thus the kernel is spanned by elements of the form $\gamma(t)\omega$, where $\gamma\in\pi_1(\mathbb C^*)$.   If $\omega_1=a^2\omega_2$, then $S_{\omega_1}=a S_{\omega_2}$ in the total space of $K^2$.  Thus above $\gamma(t)\omega$ lies the curve $\sqrt{\gamma(t)}S_\omega$, showing that if $\gamma$ is a generator for $\pi_1(\mathbb C^*)$, then $\gamma\omega$ acts by the involution $\tau$.  
\end{proof}

Notice that since the degree of $K$ is $2g-2$ (Gauss-Bonnet, e.g), $K^2$ is of degree $4g-4$.  
Let $M^{[4g-4]}$ be the configuration space of $4g-4$ distinct, unordered points on $M$.  $\pi_1(M^{[4g-4]})$, the surface braid group, has been studied in great detail, see  \cite{Bellingeri} \cite{Birman}, \cite{Scott}; we also discuss some necessary results in Section \ref{sec.surf}.   Let $\rho:H^s\to M^{[4g-4]}$ be the map which takes a quadratic differential to its zero set.  This map factors through $PH^s\hookrightarrow M^{[4g-4]}$, which is an injection by the second paragraph of this section. 

\begin{proposition}\label{prop.inj}
$\rho_*\pi_1(H^s) <\text{ker}(\pi_1(M^{[4g-4]})\to H_1(M,\mathbb{Z}))$.
\end{proposition}

\begin{proof}
This proposition will follow very quickly, once we've defined the map \\$\pi_1(M^{[4g-4]})\to H_1(M,\mathbb{Z})$.  

Define the Abel map, 
\[
A:M^{[4g-4]}\to \text{Jac}(M)\cong\text{Pic}^0(M),
\]
$A(D)=|D|\otimes K^{-2}$.  In fact, one may choose any degree $4-4g$ basepoint.  We chose $K^2$ for the simple reason that $A^{-1}(0)=PH^s$.  Therefore, 
\[
\rho_*:\pi_1(H^s)\to\ker (\pi_1(M^{[4g-4]})\to\pi_1(\text{Jac}(M))).
\]
However, $\text{Jac}(M)=H^0(M,K)^*/H_1(M,\mathbb{Z})$, thus there is a natural isomorphism, $\nu:H_1(\text{Jac}(M))\to H_1(M,\mathbb{Z})$.  Adding Hurewicz, the composition
\[
\xymatrix{
\pi_1(M^{[4g-4]})\ar[r]^{A_*}&\pi_1(\text{Jac}(M)) \ar[r]^{\text{Hurewicz}}& 
H_1(\text{Jac}(M))\ar[r]^{\nu} & H_1(M,\mathbb{Z})
}
\]
takes a braid in $M^{[4g-4]}$, and closes it to a union of homology classes of loops in $M$.  The theorem follows from the observation that $A\circ\rho=0$. 
\end{proof}

We will show in Theorem \ref{thm.iso} that for $g\geq 3$ this is an isomorphism.
\section{Surface braid groups}\label{sec.surf}
Here we collect results and definitions concerning the surface braid group.  Let $R$ be any (compact) Riemann surface.  

Let $R^n$ be the configuration space of distinct ordered subsets of $R$ of cardinality $n$.  Denote by $R^{[n]}$ the configuration space of $n$-tuples of distinct \textit{unordered} points on $R$.  $R^{[n]}$ is a quotient of $R^n$ by the symmetric group.  Let $R_m=R\setminus X$ for some $X\subset R$ of order $m$.  For any choice, $X$, the spaces $R_m$ are homoemorphic, and we intend to study only topological data.  Let $R_m^n=(R_m)^n$ and $R_m^{[n]}=(R_m)^{[n]}$.  We will also use the notation $R^{[m]+[n]}$ to denote the space of all $X,Y\subset R$ which are disjoint and of cardinality $m$ and $n$ respectively.  

The following theorem is classical:
\begin{theorem}[Fadell-Neuwirth \cite{Fadell}, \cite{Birman}]
The following is a fibration:
\[
R^n_m\to R^{n+m}\to R^m.
\]
\end{theorem}

We will use two corollaries.

\begin{corollary}\label{cor.fn1}
\[
\pi_1(R_1^{[n]})\lhd\pi_1(R^{[n]+1}).
\]
\end{corollary}

\begin{corollary}\label{cor.fn2}
\[
\pi_1(R_1^{[n]})\pi_1(R_n^1)\twoheadrightarrow \pi_1(R^{[n]+1}).
\]
\end{corollary}

We will also use part of the structure theory of surface braid groups:

\begin{theorem}\label{thm.structure}
Let $U$ be an open disk in $M$, $m\in M\setminus U$, $X\in U^{[n]}\subset M^{[n]}$.  Then
\begin{enumerate}
\item $\pi_1(U^{[n]},X)$ is generated by transpositions $\sigma_1,\ldots,\sigma_{n-1}$ of $x_j$ and $x_{j+1}$ (planar structure theory).
\item $\pi_1(M\setminus\{m,x_2,x_3,\ldots,x_n\},x_1)$ is generated by $\sigma'_\alpha$ for $\alpha$ in some basis of  $H_1(M\setminus m,\mathbb Z)$, where $x_1$ traces a path homotopic in $M$ to $\alpha$.
\item $\pi_1((M\setminus m)^{[n]},X)$ is generated by $\sigma_1,\ldots,\sigma_{n-1}$ from $(1)$, and $\sigma_\alpha$ from $(2)$.
\end{enumerate}
\end{theorem}

Here we also give a brief technical definition for transposition.
Let $e$ be an interval embedded in $M$ with endpoints $x_1$, $x_2$.  Let $X=\{x_1,\ldots,x_n\}\in M^{[n]}$, and assume also that $x_j\notin e$ for $j>2$.  If $U$ is some contractible neighborhood of $e$, $X\cap U=\{x_1,x_2\}$, then
\[
\pi_1(U^{[2]},\{x_1,x_2\})\cong \mathbb{Z}.
\]
We let $s=\{s_1,s_2\}:[0,1]\to U^{[2]}$ be either generator of this group.  Extending by trivial paths $s_j\equiv x_j$, we get an element  $\sigma\in\pi_1(M^{[n]})$.  We call the elements $\sigma_e$ and $\sigma_e^{-1}$ \textbf{transpositions} associated to $e$.  

\section{Hyperelliptic Curves}\label{sec.hyperelliptic}
Beginning in this section, we assume that $M$ is hyperelliptic.  Notationally, there is some polynomial $f$ of degree $2g+2$ with distinct roots such that $M$ is the zero set $y^2=f(x)$.  Without loss of generality, $f(0)\neq 0$.  This is a two-to-one cover $M\to \mathbb{P}^1$ ramified at the $2g+2$ roots of $f$.  Label (arbitrarily) the points in $M$ above infinity as $\infty^\pm$, and the points above $0$ as $0^\pm$.  These points will ultimately be the vertices of some graph, $\Gamma$.  $M$ is equipped with a hyperelliptic involution $\iota:(x,y)\to (x,-y)$. $\iota$ gives a decomposition of the space of quadratic differentials:
\[
H^0(M,K^2)=H^0(M,K^2)^+\oplus H^0(M,K^2)^-.
\]

The space of quadratic differentials is generated by $\dxys$, which is holomorphic with divisor $\langle\dxys\rangle=(2g-2)\langle\infty^+\rangle+(2g-2)\langle\infty^-\rangle$.  This element is fixed under the hyperelliptic involution.  

We may decompose any quadratic differential:
\[
\omega=\omega^++\omega^-=p(x)\dxys +q(x)y\dxys,
\]
such that $\iota\omega^\pm=\pm\omega^\pm$.  The degree of $p$ is at most $2g-2$ and the degree of $q$ is at most $g-3$.  The dimensions of the components of $H^0(M,K^2)$ are $h^0(K^2)^+=2g-1$, and $h^0(K^2)^-=g-2$.  In the genus $2$ case, all quadratic differentials are even.

\section{genus 2}
We treat first the special case $g=2$.  This is quite simple and unique because $H^0(M,K^2)=H^0(M,K^2)^+$.  Thus a quadratic differential is described by the polynomial $p$ of degree at most $2$.  Let $D_5^{[2]}$ be the configuration space of pairs of distinct unordered points on the five-punctured disk.

\begin{theorem}
If $M$ is a hyperelliptic curve of genus $2$, then $\pi_1(PH^s)$ is the group, $\pi_1(D_5^{[2]})$, of 2-braids on the 5-punctured disk, and $\pi_1(H^s)\cong \mathbb Z\times \pi_1(D_5^{[2]})$.
\end{theorem}

\begin{remark}
In fact, we show that $PH^s$ is merely the set of pairs of distinct points on $\mathbb{P}^1\setminus f^{-1}(0)$, which is $(\mathbb{P}^1)^{[2]}_6=D_5^{[2]}$.
\end{remark}

\begin{proof}
Notice that $H^0(M,K^2)$ consists of quadratic differentials $\omega=p(x)\dxys$, where $p$ is of degree at most two.  Such a polynomial has two zeroes in $\mathbb{P}^1$.  These lift to four zeroes in $M$, which will be unique exactly when the two zeroes are distinct and miss the ramification divisor of $M\to \mathbb{P}^1$.

Choosing some $x_0$ in the ramification divisor, $p(x_0)\neq 0$, so $\omega\to (p(x_0),\langle\omega\rangle)\in \mathbb C^*\times D_5^{[2]}$ gives a splitting of $H^s\cong \mathbb C^*\times D_5^{[2]}$, so this shows $\pi_1(H^s)\cong \mathbb Z\times \pi_1(D_5^{[2]})$.
\end{proof}

\section{Cellular decomposition of M}\label{sec.cellular}
The objective of the next two sections will be to describe explicitly the structure of $\rho_*\pi_1(H^s)$ as in Theorem \ref{thm.iso}.  We prove this theorem by applying a theorem from \cite{Copeland}:

\begin{theorem}\label{thm.poly}
If $M$ is a polyhedron (two-dimensional cell complex) of genus $g$ with $n$ faces such that no face is a neighbor of itself and no two faces share more than one edge, then $B_n^0=\ker(\pi_1(M^{[n]})\to H_1(M))$ is generated by the edge set.  Specifically, the basepoint of $M^{[n]}$ may be chosen to be a marked point in the interior of each face, and each edge may be viewed as a transposition of the marked points on the faces it separates.
\end{theorem}

We will give a cellular decomposition of $M$, and then show all generators from Theorem \ref{thm.poly} lie in $\rho_*\pi_1(H^s)$, thus proving Theorem \ref{thm.iso}.  This argument relies heavily upon $M$ being hyperelliptic.  Also, this is the only place we actually need to use that $M$ is hyperelliptic, so this is the step needed to generalize to any curve (see Section \ref{sec.generalize}).

First of all, notice that a polyhedron of genus two with four faces must necessarily have at least seven edges by Gauss-Bonnet.  However, a polyhedron with four faces conforming to the hypotheses of Theorem \ref{thm.poly} may have no more than six edges.  Thus in the case $g=2$, we would not be able to apply this theorem.

Now assume $g>2$.  We begin by constructing a triangulation of $M$ by giving a graph $\Gamma_T$ on $M$.

For simplicity we restrict our attention to the case $f(x)=x^{2g+2}-1$, from which it should be clear how to generalize.  In the closing remarks of Section \ref{sec.g3}, we discuss how one would generalize.  

On $\mathbb{P}^1$, $M$ branches over the set $\{x|0\leq x^{2g+2}\leq 1 \}$.  That is, the complement of this set in $\mathbb P^1$ lifts to a disjoint pair of sets in $M$.
Draw the set $\{x|-\infty\leq x^{2g+2}\leq 1 \}$.  This is a graph $\Gamma_T$ with vertices $0^\pm$ and $\infty^\pm$, as in Figure 1, which gives a triangulation of $M$ into $4g+4$ faces.  The interior of each triangle lies entirely in the upper or lower branch of $M$.  Triangles in the upper (e.g.) branch have vertices $(\infty^+,0^+,0^-)$.  The dual graph $\check{\Gamma}_T$ to (the fat structure on) $\Gamma_T$ is shown in Figure 2.  Its vertices are triangles in the triangulation.  In such a way, we get a triangulation for $M$ hyperelliptic of genus $g\geq 3$ (indeed, for $g\geq 1)$.  Call this the \textbf{fundamental triangulation of M}.

$\check\Gamma$ adopts a fat graph structure from $M$, and one sees that the cyclic ordering on the vertices of $\check\Gamma$ is as in Figure 2 for the upper (outer) vertices, and is reversed for the lower (inner) vertices.  

\begin{figure}[htb]
\ \ \ \ \ \ \
\includegraphics[width=.3\textwidth]{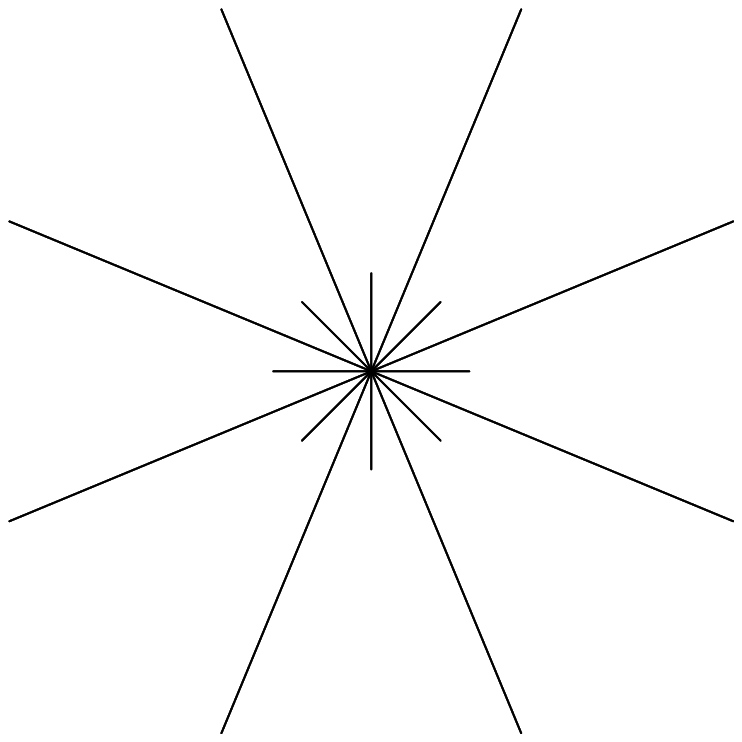}
\ \ \ \ \ \ \ \ \ \ \ \ \ \ \ \ \ \ \ \ 
\includegraphics[width=.3\textwidth]{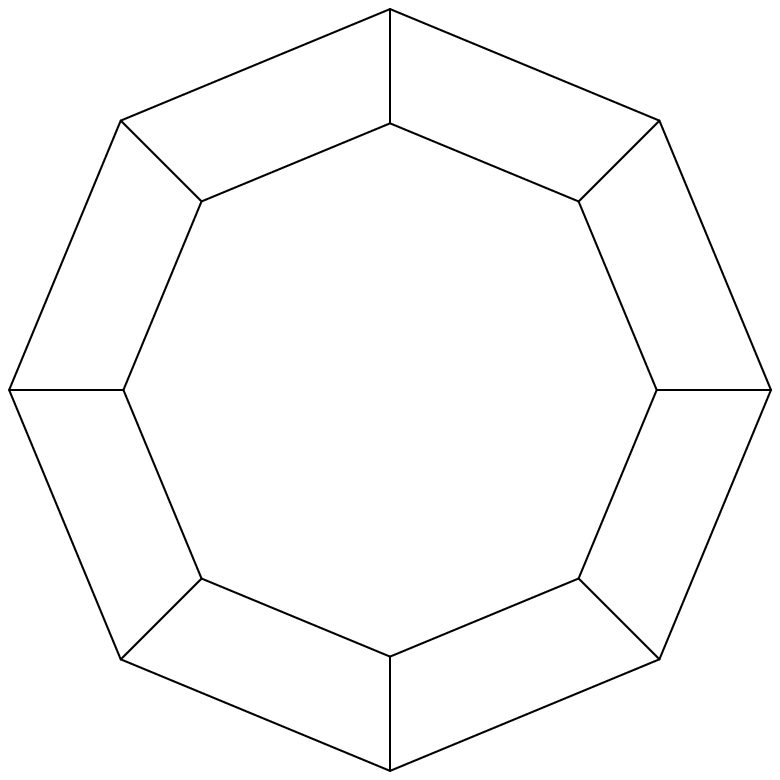}

Figure 1. Triangulation of genus 3 curve.\ \ \ \ \ \ \ \ \ Figure 2. $\check{\Gamma}_T$ for genus 3.\ \ \ \ \
\end{figure}

Notice that if two triangles are neighbors in a triangulation across an edge, $e$, then removing $e$ from the edge set replaces these two triangles with one quadrilateral.  Beginning from the fundamental triangulation which has $4g+4$ faces, we will erase 8 edges and join 16 triangles to form 8 quadrilaterals.  The new polyhedron will have $4g-4$ faces.

The operation of erasing an edge is easier to envision on $\check\Gamma$.  Erasing an edge collapses two vertices of $\check{\Gamma}_T$ along some edge.  Such an operation is admissible (with respect to the hypotheses of Theorem \ref{thm.poly})  if it creates no loops or double edges.  In the next section, we will also use the existence of some \textit{even} quadratic differential $\omega$ such that the interior of each face contains exactly one zero of $\omega$.  From this data, we may apply Theorem \ref{thm.poly}.

Figures 3, 4, and 5 illustrate which edges are erased for the genus 3, 5, and 10 surfaces.  The eight gray lines will be collpsed.

\begin{figure}[htb]
\noindent
\includegraphics[width=.28\textwidth]{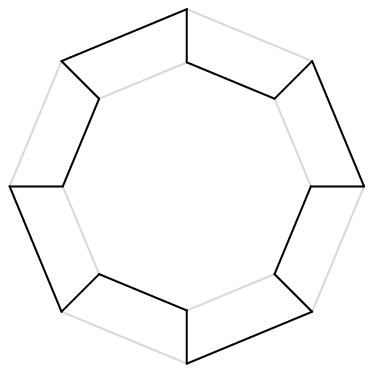}
\ \ \ \ \ \
\includegraphics[width=.28\textwidth]{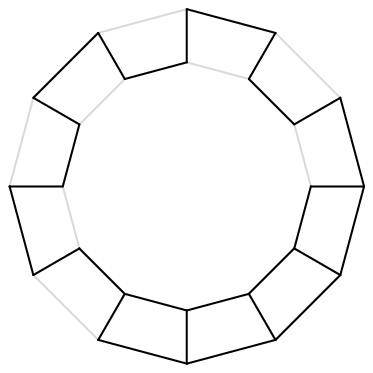}
\ \ \ \ \ \ \
\includegraphics[width=.28\textwidth]{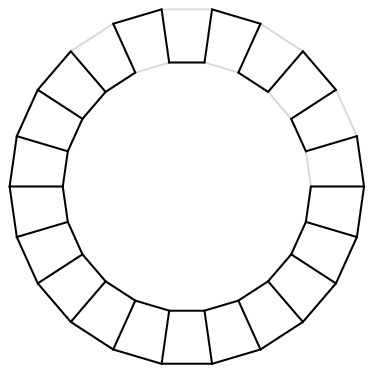}

\ \ \ Figure 3. \ \  \ \ \ \ \ \ \ \ \ \ \ \ \ \ \ \ \ \ \ \ \ \ Figure 4. \ \ \ \ \ \ \ \ \ \ \ \ \ \ \ \  \ \ \ \ \ \ \ \ \ Figure 5.\ 
\end{figure}

Explicitly, the edges of $\Gamma_T$ lying on the upper branch (as well as the lower branch) are parametrized by $re^{2\pi i(k+1/2)/(2g+2)}$, with $r>0$.  Label these edges by $k$, and likewise for the lower branch.  
On the upper branch, we remove the edges labelled by $k=1,3,5,7$ and on the lower branch, remove the edges labelled by $k=0,2,4,6$.  Let $\Gamma$ and $\check\Gamma$ be the new graphs created from $\Gamma_T$ and $\check\Gamma_T$ by these deletions.  Figures 6, 7, and 8 illustrate $\check{\Gamma}$ after applying this operation in $g=3,5,10$. 

\begin{figure}[htb]
\includegraphics[width=.28\textwidth]{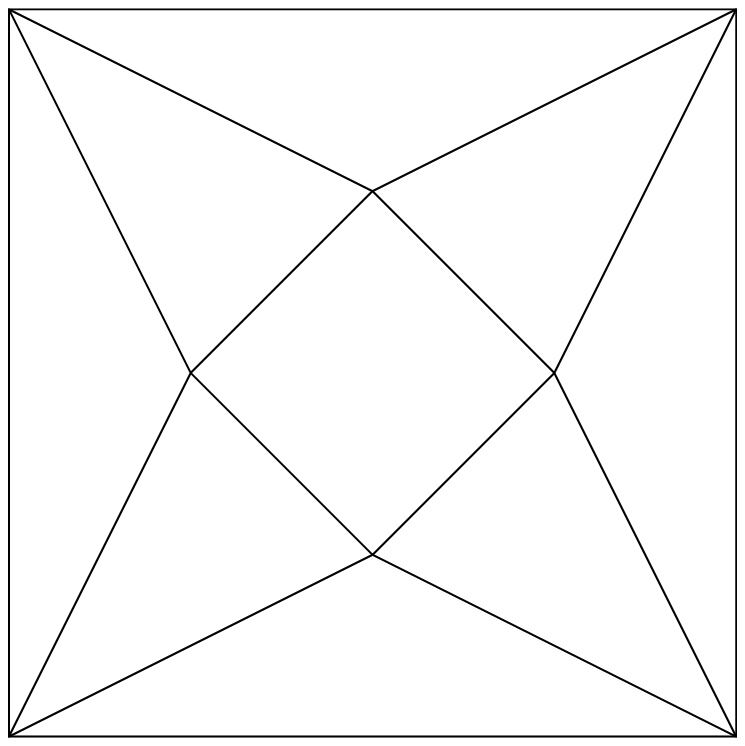}
\ \ \ \ \ \
\includegraphics[width=.28\textwidth]{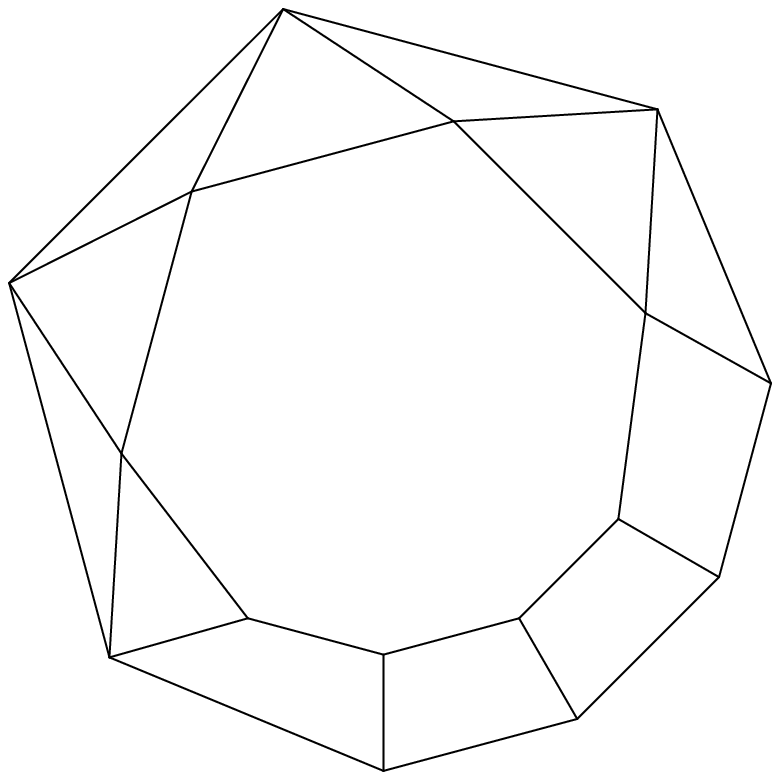}
\ \ \ \ \ \
\includegraphics[width=.28\textwidth]{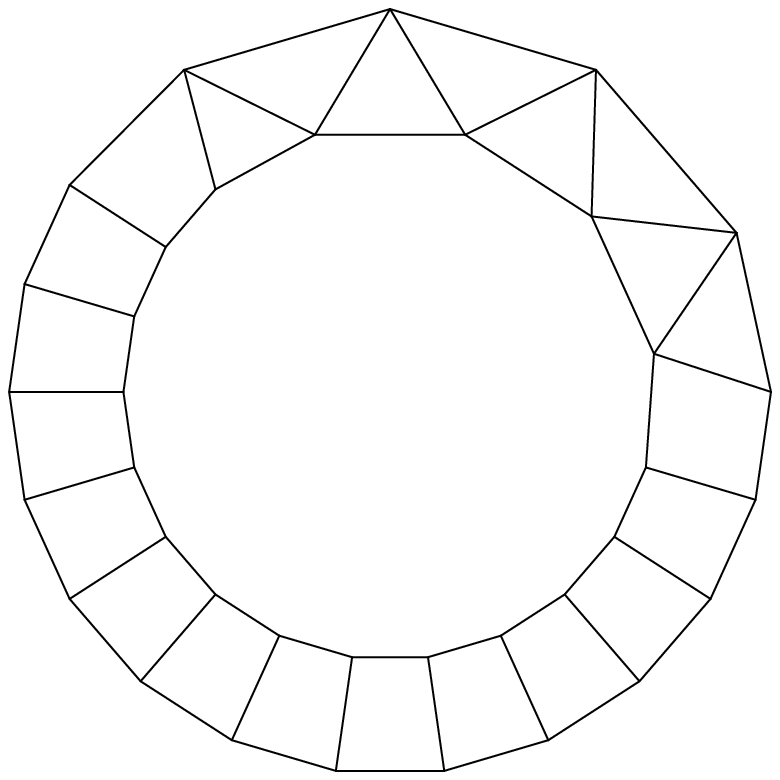}

\ \ Figure 6.\ \ \ \ \ \ \ \ \ \ \ \ \ \ \ \ \ \ \ \ \ \ \ \ Figure 7. \ \ \ \ \ \ \ \ \ \ \ \ \ \ \ \ \ \ \ \ \ \ \ Figure 8.
\end{figure}

These graphs have $4g-4$ vertices, no loops and no double edges.  Thus this gives an appropriate cellular decomposition for $M$.  

\section{genus 3 and larger}\label{sec.g3}

We now prove the structure theorem for $\pi_1(H^s)$.  Recall that in Proposition \ref{prop.inj}, we showed that $\rho_*\pi_1(H^s)$ is a subgroup of $\ker(\pi_1(M^{[n]})\to H_1(M))$.  In fact for genus greater than 2, this is an equality:

\begin{theorem}\label{thm.iso}
If $M$ is a hyperelliptic curve of genus $g>2$, then $\rho_*\pi_1(H^s)=\ker(\pi_1(M^{[n]})\to H_1(M))$.  
\end{theorem}

This theorem, coupled with Corollary \ref{cor.1tau} allows us to circumvent discussion of generators for $\pi_1(H^s)$ when studying the monodromy.  Applying Theorem \ref{thm.poly}, we get the corollary:

\begin{corollary}\label{cor.edgesgenerate}
The transpositions associated to the edge set of $\check\Gamma$ generate the image $\rho_*\pi_1(H^s)$.
\end{corollary}

This observation is in fact the method of proof for the theorem.

\begin{proof}[Proof of Theorem \ref{thm.iso}]
Let $K=\ker(\pi_1(M^{[n]})\to H_1(M))$. 
We know by Proposition \ref{prop.inj} that $\rho_*\pi_1(PH^s)$ is a subgroup of $K$.  By Section \ref{sec.cellular}, we may treat $M$ as a $(4g-4)$-hedron topologically.  We will find a basepoint  $\omega \in H^s$ such that each face of $(M,\Gamma)$ contains one zero of $\omega$.  Applying Theorem \ref{thm.poly} we need to show that all generators of $K$ are in $\rho_*\pi_1(H^s)$.  Thus we need to construct the transposition associated to each edge of $M$ (technically one of two possible transpositions).

Recall that $H^0(M,K^2)=H^0(M,K^2)^+\oplus H^0(M,K^2)^-$, and $\text{dim}( H^0(M,K^2)^-)>0$ for $g\geq 3$.  Let $\zeta=e^{2\pi i/2g+2}$.
  Notice that $x\mapsto x^{2g+2}-2^{2g+2}$ has one root in each triangle of the fundamental triangulation of $M$ (recall Figure 1).  Consider the quadratic differential 
\[
\omega=(x-2 \zeta^2)(x-2\zeta^4)(x-2\zeta^6)(x-2\zeta^8)\prod_{9\leq j\leq 2g+2}  (x-2\zeta^j)\dxys.
\]
This has $4g-4$ simple zeroes, and each face of $M$ contains exactly one zero. Figure 9 shows the zeroes of $\omega$ in the fundamental triangulation of $M$ for $g=5$.  The rays are labelled by $u$, $l$, or $u/l$, for whether they represent edges in the upper or lower branch of $M$ or both.

\begin{figure}[htb]
\includegraphics[width=.3\textwidth]{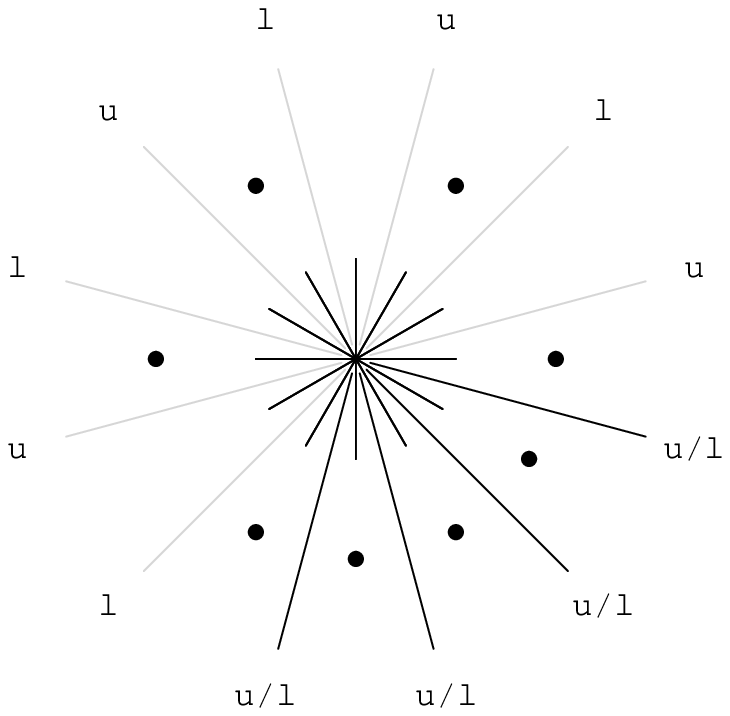}

Figure 9. 
\end{figure}

Finally, we must show that transpositions across the edges are contained in $\pi_1(H^s)$.  The remainder of this section is a string of lemmas dealing with this issue.
\end{proof}

The reader should be warned that the remainder of this section is quite computational.  One might prefer to skip to the last paragraph of this section on first reading.

Let $\kappa=\sqrt{2^{2g+2}-1}>0$.  The roots of $\omega$ are all of the form $(x,y)=(2\zeta^j,\pm\kappa)$.  For the remainder of this section, we fix $\omega$, $\kappa$, $\zeta$ as above.  As our convention, we will assume that the ``upper branch'' of $M$ contains the points $(2\zeta^j,(-1)^j\kappa)$.

\begin{lemma}
Fix $j$, $1\leq j\leq 2g+1$.  The pair of points $(x,y)=(2\zeta^j,\pm\kappa)$ are neighbors across some edge $e_\alpha$.  Their transposition is an element of $\rho_*\pi_1(H^s)$.
\end{lemma}

\begin{proof}
The two faces containing $x=\zeta^j$ are parametrized by complex numbers, $x$, with argmument $\frac{j-1/2}{2g+2}<\theta<\frac{j+1/2}{2g+2}$.  Call their union $D$.   One may use $y$ as a coordinate on the interior of $D$.  In this coordinate, $D=\{y|y^2\notin \mathbb{R}_{<-1}\}$, as in Figure 10.

\begin{figure}[htb]

\includegraphics[width=.3\textwidth]{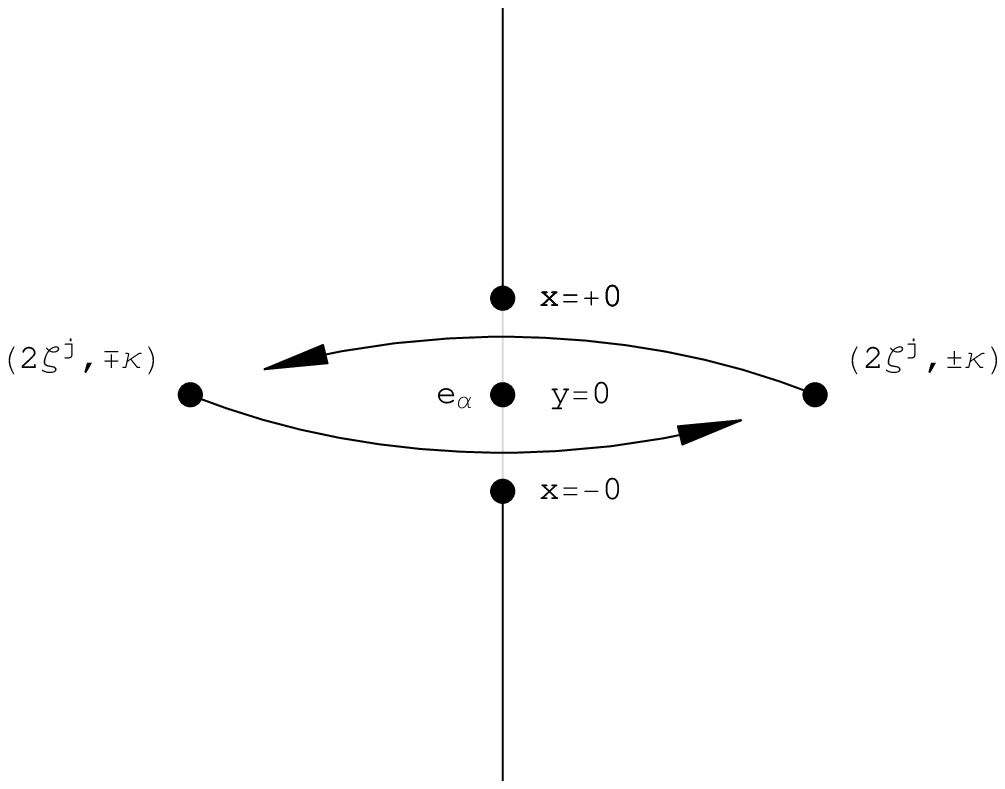}

Figure 10.
\end{figure}

Let $\omega(t)=\frac{\omega}{x-2\zeta^j}(x-\alpha(t))$, where $x=\alpha(t)$ is an arc in $D$ connecting $y=+\kappa$ to $y=-\kappa$.  $[\omega(t)]\in\rho_*\pi(PH^s)$ is the transposition.  As $\alpha$ changes, the points trace a pair of paths symmetric about $y=0$.  One may choose $\alpha$ explicitly as:
\[
\alpha(t)=2\zeta^j(\frac{1-\epsilon}{2}+\frac{1+\epsilon}{2}\cos(t)+i\sin(t))
\]
for any small enough $\epsilon$ such that the range of $\alpha$ is in $D$.  This gives the curve from Figure 10.
\end{proof}

\begin{lemma}\label{lem.stupid}~
\newline
\begin{enumerate}
\item Let $\lambda=\pm \kappa$.   If $a=(2\zeta^j,\lambda)$ and $b=(2\zeta^{j+1},-\lambda)$ are roots of $\omega$, then the transposition of $a$ and $b$ is in $\rho_*\pi_1(H^s)$.
 \item Let $\lambda=\pm \kappa$.   If $a=(2\zeta^j,\lambda)$ and $b=(2\zeta^{j+2},\lambda)$ are roots of $\omega$ and $(2\zeta^{j+1},\lambda)$ is not a root, then the transposition of $a$ and $b$ is in $\rho_*\pi_1(H^s)$.
\end{enumerate}
\end{lemma}

\begin{proof}
We prove (1).  (2) is identical.  The points $a$ and $b$ are neighbors in $M$.  Examine all curves:

\[
\omega_\epsilon(t)=\frac{\omega}{x-2\zeta^j}(x-2+\frac12\zeta^{j+\sin(\pi t)}\sin(2\pi t))+\epsilon t(1-t)y\dxys.
\]
for $\epsilon\in\mathbb C$.  For small enough $\epsilon$, the trajectories of all points other than $2\zeta^j$, and $2\zeta^{j+1}$ are constrained to their faces.  If we choose $\epsilon$ such that 
\[
\left.\frac{~~~\frac{\epsilon y}4\dxys~~~}{~~~\frac{\omega}{x-2\zeta^{j}}~~~}\right|_{(x,y)=(2\zeta^{j+1},-\lambda)}
\]
is a small positive multiple of $i\zeta^{j+1}$, then $[\omega(t)]$ is the desired transposition.

\end{proof}

\begin{lemma}
Assume that $a=(2\zeta^j,\kappa)$ and $b=(2\zeta^{j+2},-\kappa)$ are roots of $\omega$ and $(2\zeta^{j+1},\pm\kappa)$ are not a roots.  Furthermore, let $e_\alpha$ be the edge of $\Gamma$ containing the branch point $\zeta^{j+1}$.  Then $a$ and $b$ are neighbors across $e_\alpha$, and the transposition associated to $e_\alpha$ is in $\rho_*\pi_1(H^s)$.
\end{lemma}

\begin{proof}
If we let $\omega_\alpha(t)=\frac{\omega}{x-2\zeta^j}(x-\alpha(t))$, where $\alpha$ is a curve turning clockwise once around $\zeta^{j+1}$, then conjugating the curve from Lemma \ref{lem.stupid} (2) by $[\omega_\alpha]$, we get the desired result.
\end{proof}

All edges of $\check{\Gamma}$ are of the types described in the above lemmas, so we are done.

All of the constructions from this section pass to general hyperelliptic curves as follows:

Let $y^2=f(x)$ be any hyperelliptic curve of genus $g>2$.  Without loss of generality, $\text{deg}(f)=2g+2$, and $f(0)\neq 0$.  Let $\omega\in H^s$.  Choosing non-intersecting paths from $0$ to $\langle\omega\rangle$, and $0$ to $\infty$, which alternate in the appropriate way, we may again construct a fundamental triangulation of $M$, and consequently construct $\Gamma$ on $M$.  The constructions of the transpositions depended only on curves in $(\mathbb{P}^1)^{2g-2}$ with small deformations proportional to $y\dxys$.

\section{Homology of $S$ via $\check\Gamma$}\label{sec.ho}
In this section, we intend to study the cellular homology of $S$ via $\check\Gamma$.  We focus also on the $\tau=-1$ part of the homology, as it is related to the homology of $\text{Prym}(S,M)$.  Let $X$ be the ramification divisor of $S\to M$.  As such, $X$ is a fixed point set for $\tau$ and each face of $\Gamma$ contains exactly one point of $X$.  Choose some realization of $\check\Gamma$ in $M$ such that the vertices of $\check\Gamma$ are $X$, each edge of $\check\Gamma$ crosses $\Gamma$ only once, and this crossing is through its corresponding edge of $\Gamma$.

We continue to use $\omega$ as a basepoint for $H^s$ and let $S=\sqrt\omega$. Recall that on $M$ we have the graph $\Gamma$ which makes $M$ a polyhedron.  The ramification points of $S\to M$ lie away from $\Gamma$, so that we may lift to $\bar\Gamma$, a graph on $S$.  Every face, $f$, of $\Gamma$ contains exactly one ramification point.  Thus if $f$ is a triangle, then its lift, $\bar f$, is a hexagon.  Each pair of opposite edges of $\bar f$ lies above a single edge of $f$.  Likewise quadrilaterals lift to octagons with opposite edges identified under $S\to M$.  Opposite vertices are also identified.  

Thus $\bar\Gamma$ has the same face set as $\Gamma$ with the edge and vertex sets doubled.   As a check, the Euler characteristic is $2(4)-2(6g-2)+(4g-4)=8-8g=2-2(4g-3)$.  This agrees with the Riemann-Hurwitz formula, which says that the genus of $S$ is $\tilde g=4g-3$.  

We may also study $\check{\bar\Gamma}$, the dual to $\bar\Gamma$ on $S$, or equivalently, the lift of $\check\Gamma$. Its  1-skeleton has the same vertex set as $\check\Gamma$ and all edges doubled.  We use the cell structure of $\check{\bar\Gamma}$ as a fat graph, or cellular decomposition of $S$, as our cell complex $C_\bullet$ for $S$.
\[
C_\bullet=(
0\to\mathbb{Z}\bar F\to\mathbb{Z}\bar E\to\mathbb{Z}X\to 0 
).
\]

Notice that the differential $\partial$ of $C_\bullet$ is the standard boundary operator.  Since $\tau$ is a cellular diffeomorphism, $\partial\tau=\tau\partial$.

\section{The graph $\Gamma$}

In this section, we restrict our attention to the case $g=10$ when drawing pictures, though all arguments go through for any  $g>2$.  Herein, we draw the graph $\Gamma$.  In Section \ref{sec.phi}, we will use this to construct $\bar\Gamma$.

Recall that the fat graph structure on $\check\Gamma_T$ agrees with the planar representation in Figure 2 for the upper vertices ,and is opposite for the lower vertices.  Once one contracts the appropriate edges, one finds the cyclic orientations for the vertices in Figures 6, 7, and 8 are still clockwise for the upper vertices and counter-clockwise for the lower vertices.  From this, one can easily draw the boundaries of the faces by following the directed cycles which pass through vertices by entering along one half-edge and leaving along the ``next'' half edge.  There are four such cycles, one for each face of $\check\Gamma$. 

The graph $\Gamma$ has four vertices labelled $0^\pm$ and $\infty^\pm$.  There are $4g-4$ edges with multiplicity.  To find the neighbors of, for example, $0^\pm$, one need only draw a small loop $\gamma(t)=\frac{1}{2}e^{it}$, $t\in [0,\pi]$, around $0$ and enumerate the edges crossed by the two curves $c^{-1}\gamma\subset S$.  Alternatively, one could find two consecutive edges of $0^+$ in $\check\Gamma$, and find which unique cycle from the previous paragraph follows these edges.  Doing this, one represent $\Gamma$ as in Figure 11.  In this figure, we've restricted to $g=10$, and labelled each multiple edge of $\Gamma$ by its dual subgraph in $\check\Gamma$.  We should point out that in the special case of $g=3$, the only neighbor of $0^-$ is  $0^+$ as in Figure 12. 

\begin{figure}[htb]
\noindent
\includegraphics[width=1.0\textwidth]{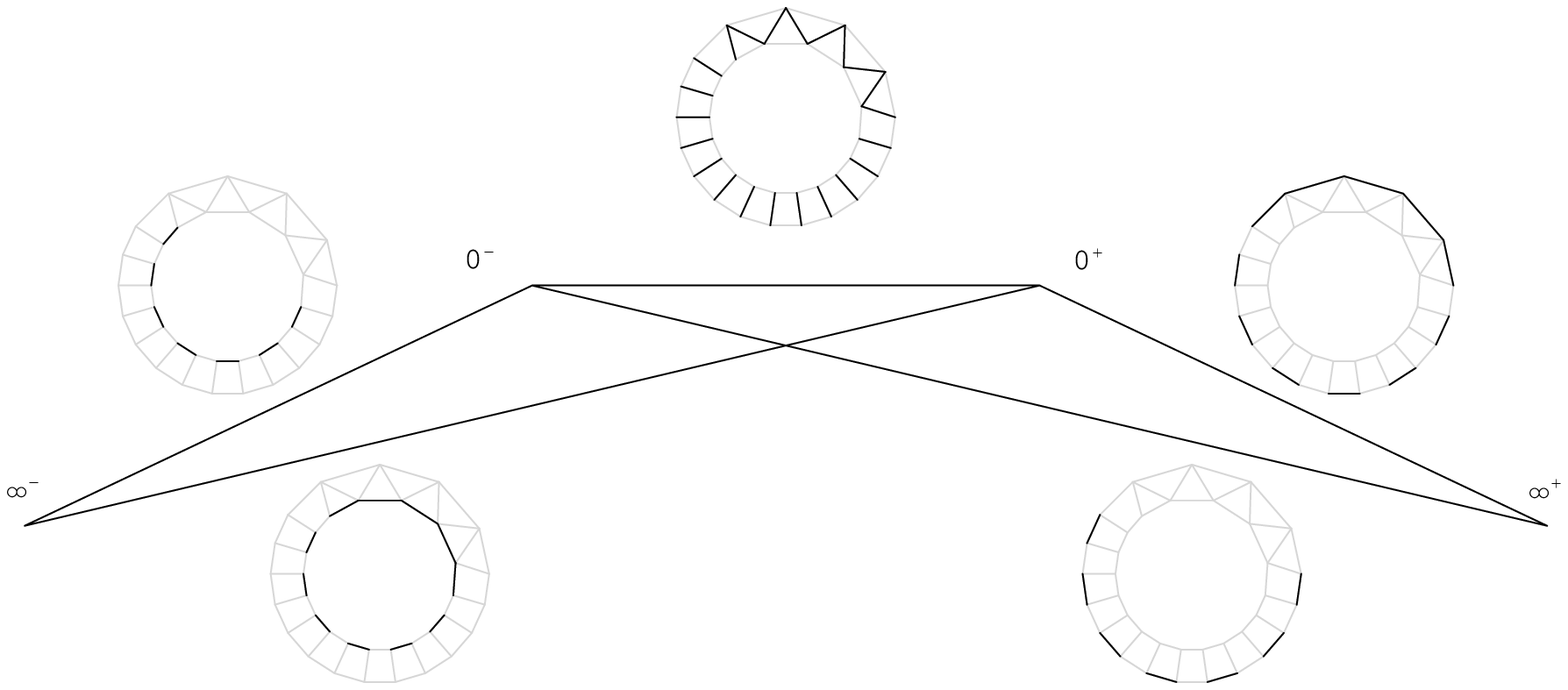}

\ \ \ \ \ \ \ \ \ Figure 11. $\Gamma$ for $g=10$.
\end{figure}

\begin{figure}[htb]
\noindent
\includegraphics[width=1.0\textwidth]{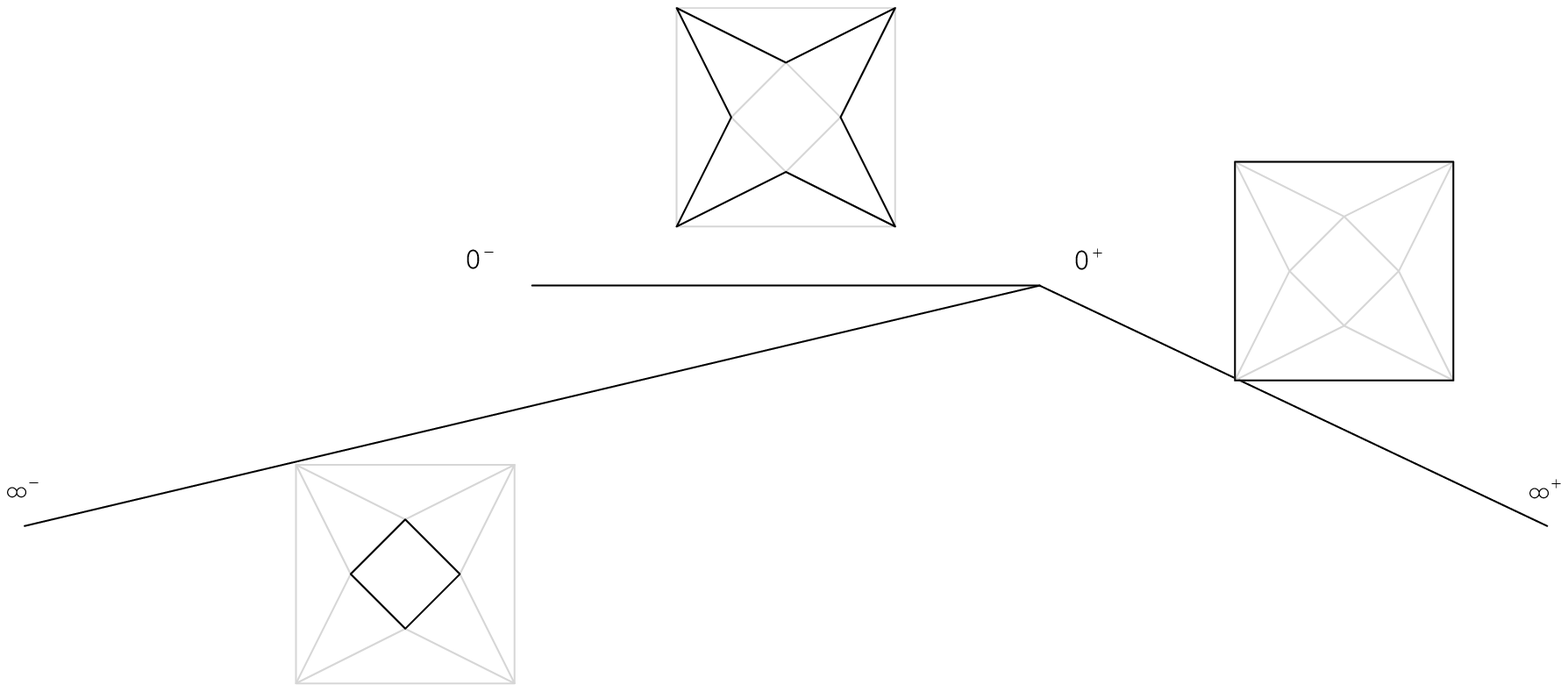}

\ \ \ \ \ \ \ \ \ Figure 12. $\Gamma$ for $g=3$.
\end{figure}

\section{orientation}\label{sec.orientation}

Since $\partial(0^++0^-+\infty^++\infty^-)=0$, there exists a unique orientation for each edge such that in $\check\Gamma$, $-\partial \infty^\pm$ and $\partial 0^-$ all are positive sums of edges.  Up to a global choice of sign, this orientation is as in Figure 13.  Since $\partial^2=0$ the consecutive segments $\partial f$ should be joined head to tail.  

\begin{figure}[htb]
\noindent
\includegraphics[width=.3\textwidth]{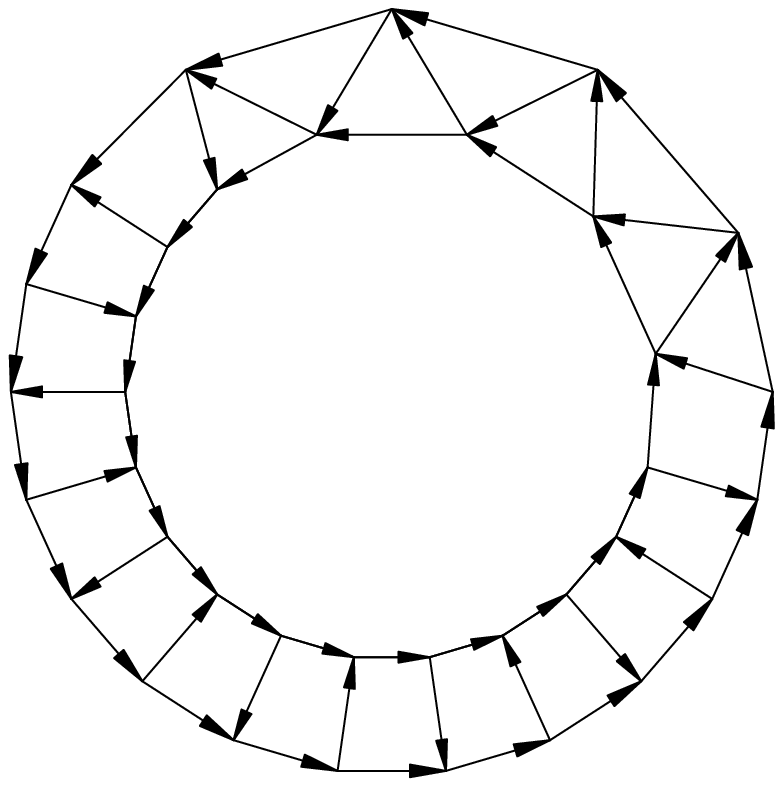}

\ \ \ \ \ \ \ \ \ Figure 13. Orienting $\Gamma$ for $g=10$.
\end{figure}

We orient the edges of $\Gamma$ so that if $e\in E$ corresponds to $e_\Gamma$ as an edge of $\Gamma$, and $e_{\check \Gamma}$ as an edge of $\check\Gamma$ (geometric realizations of $e$), then the intersection pairing on $M$ gives:
\[
e_{\check\Gamma}\cdot e_\Gamma=+1.
\]
Physically, this means that edges of $\Gamma$ travel from $0$ to $\infty$, and ``zig-zag'' between $0^+$ and $0^-$.

\section{The Prym Variety}\label{sec.prym}
In this section, we set out to convert the Prym variety into entirely combinatorial data, ultimately proving Theorem \ref{thm.prymef}.  
A basis for $H_1(S,\mathbb{Z})$ extends to an $\mathbb{R}$-basis for $H^0(S,K)^*$, as it is a full-rank sublattice.  Since we are interested in topological data (monodromy), we may ignore holomorphic structure and use the identification
\[
\text{Jac}(S)=H_1(S,\mathbb{Z})\otimes_\mathbb{Z}\mathbb{R}/H_1(S,\mathbb{Z}).
\]

Recall that $\check\Gamma$ is embedded in $M$, and is the combinatorial data $(X,E,F)$ with vertices $X$, the ramification divisor.  We lift the orientation of the previous section from $\check\Gamma$ to $\cbG$.  Any edge, $e$, of $\check\Gamma$ lifts as a set to a curve $\bar e=c^{-1}(e)$ in $S$ such for either choice of orientation on $\bar e$, $\tau\bar e=-\bar e$.  We know also that $\bar e$ is a signed sum of elements of $\cbG$.


Again pick some edge $e\in E$.  This lifts as an oriented curve to some $e_1+e_2\in \mathbb{Z}\bar E$.  $\pm(e_1-e_2)$ are the loops which lie (as sets) over $e$ and on which $\tau$ acts by $-1$.

Let $\psi:\mathbb{R}E\to\mathbb{R}\bar E$ be any map such that for $e\in E$, $\psi(e)=\pm(e_1-e_2)$, where $c(e_1)=c(e_2)=e$, $\tau e_1=e_2$.  There are $2^{|G|}$ such choices of maps.  Any such $\psi$ gives a map $\psi:\mathbb{R}E\to(1-\tau)\mathbb{R}\bar E=(\mathbb{R}\bar E)^{-}$.

Also choose a map $\Psi:\mathbb{R}F\to\mathbb{R}\bar E$ by choosing $f_1,f_2\in \bar F$ above each $f\in F$ with $\tau f_1=f_2$, and letting $\Psi(f)=\partial(f_1-f_2)$.  We will see that $\Psi\neq\psi\circ\partial$, and that $\Psi$ actually encodes data relevant to the double covering map.  Since $\psi$ is an injection, and the image of $\Psi$ lies in the image of $\psi$, we will sometimes use $\psi$ to denote $\Psi^{-1}\psi:\mathbb R F\to \mathbb RE$.

In the next section, we will construct a specific such pair $(\psi,\Psi)$, but for now we treat general $(\psi,\Psi)$.  The following theorem uses the simplicial structure of $S$ to construct the $\tau=-1$ part of $\text{Jac}(S)$.  The point of the theorem is that the Jacobian is spanned by $\bar E$.  Inside this set is $\psi E$ which is generated by elements of the form $e_1-\tau e_1$.  In the quotient $\tau=-1$, The element $\frac12(e_1-\tau e_1)=e_1$ is an honest edge, and the point of the theorem is that these elements span (over $\mathbb R$) the Prym variety.

\begin{theorem}\label{thm.prymef}
If $E$ and $F$ are the edge and face sets of $\check\Gamma$, and $(\psi,\Psi)$ is a map as above, then topologically, $(\psi,\Psi)$ induces a homeomorphism:
\[
\text{Prym}(S,M)\cong\mathbb{R}E/(\mathbb{R} F+{\scriptstyle \frac{1}{2}}\mathbb{Z}E).
\]
\end{theorem}

\begin{remark}\label{rem.prym}
One employs the isomorphism as follows.   $(\psi,\Psi):\mathbb{Z}E\times \mathbb{Z}F\to \mathbb{Z} \bar E$, which descends to an isomorphism under the identification 
\[
\text{Prym}\cong(\mathbb{R}\bar{E}/\mathbb{R}\bar F +\mathbb{Z}\bar E)^-
\]
\end{remark}

\begin{proof}
Let $K=\ker(\mathbb{R}\bar E\to\mathbb{R}X)$.
 First we show the identification of Remark \ref{rem.prym}.
\begin{align}
\text{Prym}&\cong \text{Jac}(S)^-\notag\\
&\cong [H_1(S,\mathbb{Z})\otimes\mathbb{R}/ H_1(S,\mathbb{Z})]^-\notag\\
&\cong[(\ker(\mathbb{R}\bar E\to\mathbb{R}X)/\partial\mathbb{R}\bar F)\Big/(\ker(\mathbb{Z}\bar E\to\mathbb{Z}X)/\partial\mathbb{Z}\bar F)]^-\notag\\
&\cong[(K/\partial\mathbb{R}\bar F)\Big/(\mathbb{Z}\bar E \cap K/\partial\mathbb{Z}\bar F)]^-\notag\\
&\cong[K/(\partial\mathbb{R}\bar F+\mathbb{Z}\bar E \cap K)]^-\notag
\end{align}
All of these maps descend from the identity map $\mathbb{R}\bar E\to\mathbb{R}\bar E$, thus are natural.  Now assume that $e\in \mathbb{R}\bar E$ such that $e+\tau e\in \partial\mathbb{R}\bar F+\mathbb{Z}\bar E$.  Then 
\[
(1+\tau) \partial e=\partial(1+\tau)e\in(\partial^2\mathbb{R}\bar F+\partial \mathbb{Z} \bar E)^+\subset \mathbb{R}X^+=0.
\]
Thus $\partial e\in (\mathbb{R}X)^+=0$ by Section \ref{sec.ho}, and $e\in K$.  Therefore, 
\[
[K/(\partial\mathbb{R}\bar F+\mathbb{Z}\bar E \cap K)]^-\cong[\mathbb{R}\bar E/(\partial\mathbb{R}\bar F+\mathbb{Z}\bar E \cap K)]^-
\]

There is a map
\[
[\mathbb{R}\bar E/(\partial\mathbb{R}\bar F+\mathbb{Z}\bar E \cap K)]^-\surj[\mathbb{R}\bar E/(\partial\mathbb{R}\bar F+\mathbb{Z}\bar E)]^-.
\]
The identification will be complete if we show that this map is injective.  Assume that $e\in\partial\mathbb{R}\bar F+\mathbb{Z}\bar E$, and that $e+\tau e\in \partial\mathbb{R}\bar F+\mathbb{Z}\bar E \cap K=(\partial\mathbb{R}\bar F+\mathbb{Z}\bar E) \cap K$.   Thus $e+\tau e\in K$, and $e-\tau e\in \mathbb{R}E^-\subset K$.  So $e\in K$.  Thus $e\equiv 0 \text{~mod~}(\partial\mathbb{R}\bar F+\mathbb{Z}\bar E) \cap K$, and the map is injective.

The next step of the proof is to show that the map $(\psi,\Psi):\mathbb{Z}E\times \mathbb{Z}F\to \mathbb{Z} \bar E$ induces an isomorphism:
\[
[\mathbb{R}\bar E/(\partial\mathbb{R}\bar F+\mathbb{Z}\bar E)]^-\cong\mathbb{R}E/(\mathbb{R}\partial F+{\scriptstyle \frac{1}{2}}\mathbb{Z}E).
\]
First of all, the map is well-defined, since $\partial F+{ \frac{1}{2}}\mathbb{Z}E\subset \text{ker}(1+\tau)$.  It is injective, because $\mathbb{R}E\cap(1+\tau)e(\mathbb{R}\bar F+\mathbb{Z}\bar E)\subset \mathbb{R}F+ \frac{1}{2}\mathbb{Z}E$.  If $e\in \mathbb{R}\bar E$, then $\frac{1}{2}(1-\tau)e\in\mathbb{R}E$, and since  $\frac{1}{2}(1+\tau)e\in\mathbb{R}\bar E^+$, $e\equiv \frac{1}{2}(1-\tau)e \text{~mod~}\mathbb{R}\bar F+\mathbb{Z}\bar E$.  Therefore it is surjective as well.
\end{proof}

The identification is complete.  However, it depends on the choice of $\psi$ (but not $\Psi$).  In the next section, we make specific choices for these maps, in order to perform explicit computation.

\section{The nature of $S\to M$}

In this section we fix $\omega$ as in Section \ref{sec.g3}:
\begin{align}
\omega&=(x-2 \zeta^2)(x-2\zeta^4)(x-2\zeta^6)(x-2\zeta^8)\prod_{9\leq j\leq 2g+2}  (x-2\zeta^j)\dxys\notag\\
      &=\frac{(x^{2g+2}-2^{2g+2})}{(x-2\zeta)(x-2\zeta^3)(x-2\zeta^5)(x-2\zeta^7)}\dxys\notag\\
      &=p(x)\dxys,\notag
\end{align}
and study what happens when we lift cycles of $M$ to $S$.  
Specifically, there are many double covers of $M$ ramified at the points, $X$, which can be distinguished by which closed loops lift to connected loops.  For instance, a double cover of a torus ramified at two points is uniquely determined by a branch cut between the two points.  However, the same branch points may give different covers, and there are in fact four distinct double covers of a torus with two given ramification points.  Thus we need to do a little work to see which double cover $S\to M$ is given by $\omega$, which is given by Theorem \ref{thm.separated}.

\begin{lemma}\label{lem.sqrtp1}
If $|x|<1$, then $p(x)=-r\zeta^{-16}e^{i\theta}$ for $r>1$ and $-\frac{5\pi}{6}\leq\theta\leq\frac{5\pi}{6}$.  
\end{lemma}

\begin{proof}
Clearly, $r>(2^{2g+2}-1)/2^4>1$.  Next one writes 
\[
p(x)=(x^{2g+2}-2^{2g+2})\cdot\frac{1}{(x-2\zeta)}\cdot\frac{1}{(x-2\zeta^3)}\cdot\frac{1}{(x-2\zeta^5)}\cdot\frac{1}{(x-2\zeta^7)}
\]
and observes, for example, that the argument of $\frac{1}{(x-2\zeta)}$ is in the interval 
\[
[\pi-\log{\zeta}-\frac{\pi}{6},\pi-\log{\zeta}+\frac{\pi}{6}]
\]
\end{proof}

\begin{lemma}\label{lem.sqrtp2}
There exists some $\epsilon>0$ such that on $R_\epsilon=\{(x,y)\in M\Big||x|<1+\epsilon\}$, $\sqrt{p}$ is a two-valued function with image having two components.
\end{lemma}

\begin{proof}
By continuity and Lemma \ref{lem.sqrtp1}, for small $\epsilon$, $p(R_\epsilon)$ is a connected set which is never real negative.  Thus $\sqrt{p(R_\epsilon)}$ is a pair of disconnected regions.
\end{proof}

\begin{theorem}\label{thm.separated}
Let $R=R_{\epsilon}$ as above.  Let $R_\epsilon^c=M\setminus R_\epsilon$.  Then $S|_{R_\epsilon}\to R_\epsilon$ is a disconnected double cover, and $S$ can be constructed by gluing the unique double cover of the pair of disks $S|_{R_\epsilon^c}\to R_\epsilon^c$ appropriately to $S|_{R_\epsilon}$.
\end{theorem}

Effectively what this theorem says is that as a four-fold covering $S\to M\to \mathbb P^1$, the first map is ramified outside $R_{\epsilon}$ and the second is ramified inside $R_{\epsilon}$.  

\begin{proof}
On $R_\epsilon$, one may choose a well-defined, holomorphic function $\sqrt{p}$.  Then $\sqrt{p}\dxy$ is a holomorphic one-form on $R_\epsilon$.  The graph of this one-form is one of the components of $S|_{R_\epsilon}$.  
\end{proof}

\section{The map $(\psi,\Psi)$}\label{sec.phi}

We set out in this section to construct the maps $\psi$ and $\Psi$ for general $g>2$.  First we will graph $\bar\Gamma$ as in Figures 11 and 12.  Recall the embedding of $\Gamma$ and $\omega^{-1}(0)$ in $M$ from Figure 14.  To these markings on $M$, we add the set 
\[
B_c=\{e^{it/2g+2}|t\in [0,2]\cup[4,6]\cup\bigcup_{4\leq j\leq g}[2j,2j+1]\},
\]
as in Figure 14.

\begin{figure}[htb]
\noindent
\ \ \ \ \ \ \includegraphics[width=.3\textwidth]{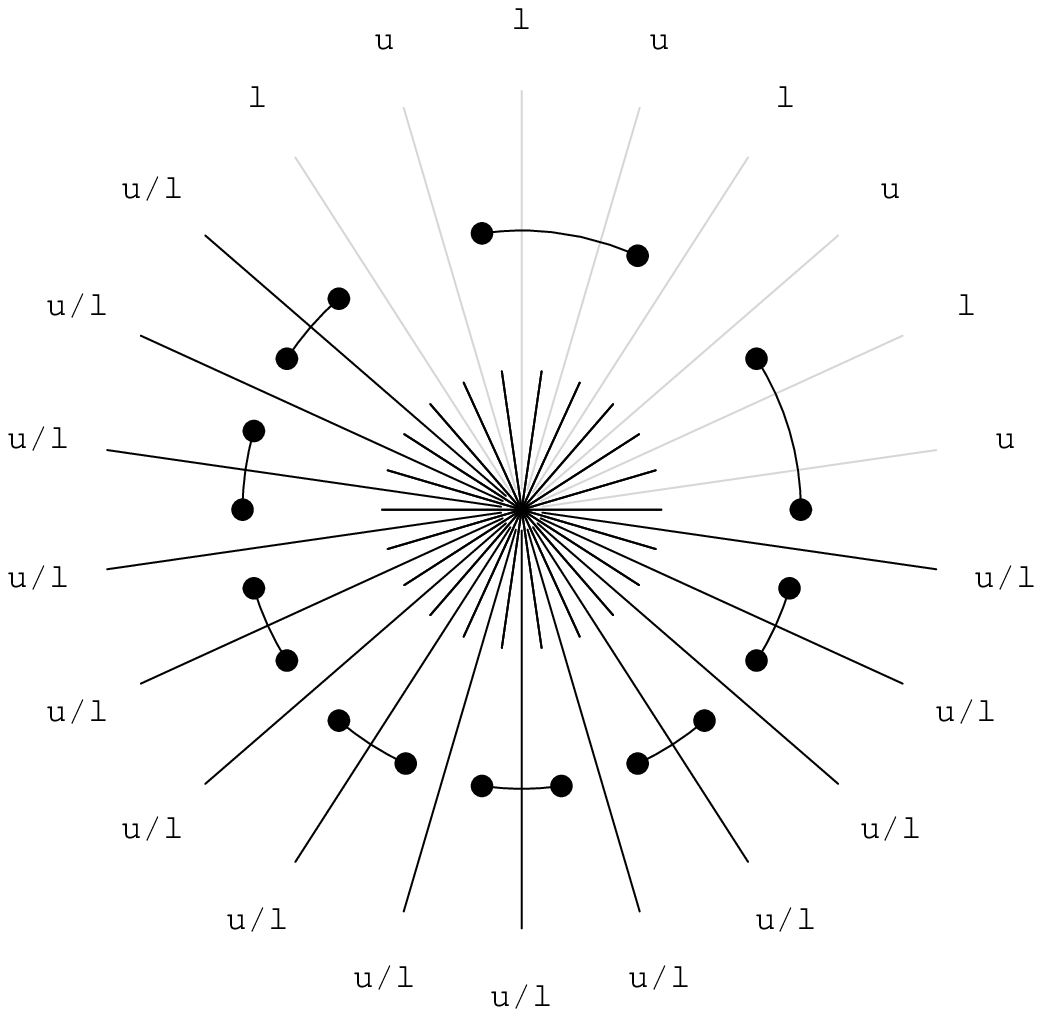}

\ \ \ \ \ Figure 14. $B_c\subset M$ for $g=10$.
\end{figure}

$B_c$ is a set of $2g-2$ arcs which may be used as branch cuts for the map $c:S\to M$.  Thus $c^{-1}(M\setminus B_c)$ is a disjoint union, $S_1\coprod S_2$, with $\tau S_1=S_2$ by Theorem \ref{thm.separated}.  Let $0^\pm_j,\infty^\pm_j$ be in $S_j$ over $0^\pm,\infty^\pm$.  Now let $\Psi(f)=\partial(f_1-f_2)$  define $\Psi$ on $F$ and extend linearly.  Clearly, if $c\bar e=e$ and $e$ connects $f$ to $f'$, then $\bar e$ connects $f_j$ to $f_k'$ with $j=k$ if $e$ does not cross $B_c$ and $j\neq k$ if $e$ does cross $B_c$.  From this observation, we may costruct $\bar\Gamma$ from $\Gamma$.  We first highlight in $\check\Gamma$ the edges that cross $B_c$ in Figure 15.  

We observe that $\Psi\neq \psi\circ\partial$ precisely due to the set $B_c$.  Working again by example, let $\bar e$ be some edge separating $\infty^+_1$ from $0^+_2$.  Then $\psi\partial(\infty^++0^+)$ is supported away from $e$.  However, since $e\in B_c$, $0_1^+$ and $\infty^+_1$ are not neighbors, $\Psi(\infty^+0^+)$ contains some component, $\pm2e$.

In our graph of $\bar\Gamma$, we will label exactly one element of each pair $\{e,\tau e\}$, where it is understood that $\tau$ gives a symmetry across the horizontal midline.  
We also choose orientations on the cycles $\psi(e)$ so that $-\Psi(\infty^\pm)$ are positive sums, and so that the cycles corresponding to edges separating $0^\pm$ appear with positive coefficient in $\Psi(0^-)$.  This is as consistent with the conventions of Section \ref{sec.orientation} as possible.

As an aside, notice that each vertex of $\check\Gamma$ (ramification point) meets $B_c$ at exactly one edge.  This is to be expected for a double cover, and will play an implicit role in Proposition \ref{prop.intnos}

\begin{figure}[htb]
\noindent
\ \ \ \ \ \ \includegraphics[width=.3\textwidth]{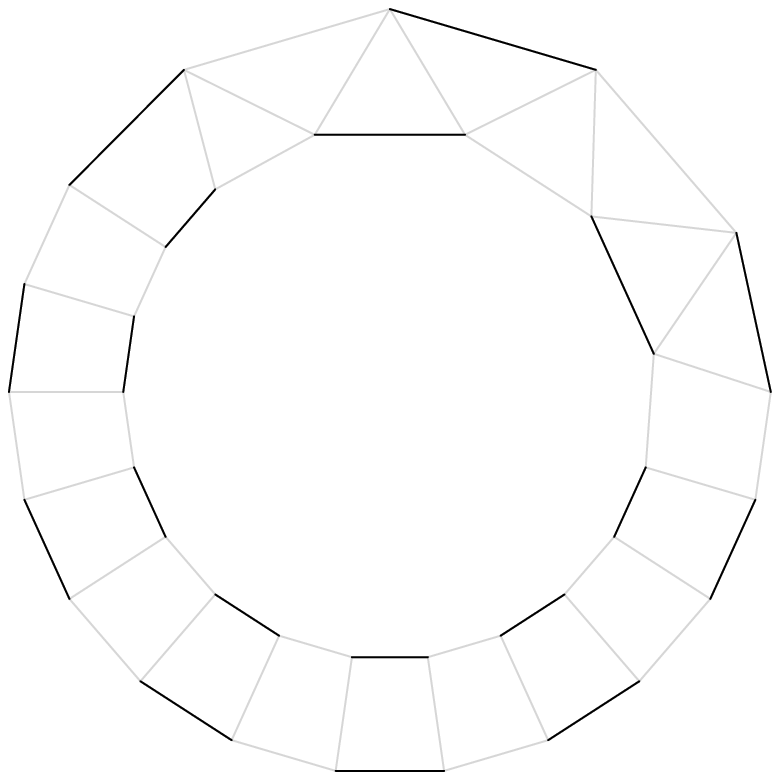}

\ \ \ \ \ Figure 15. $B_c\subset \check\Gamma$ for $g=10$.
\end{figure}

In our graph of $\bar\Gamma$, which for genus 10 is Figure 16, we will label exactly one element of each pair $\{\bar e,\tau \bar e\}$, where it is understood that $\tau$ gives a symmetry across the horizontal midline.  If for some edge $e\in E$, we have drawn $e_1$ with $ce_1=e$, then we let $e_2=\tau e_1$, and define $\psi e=e_1-e_2$.  This extends to give the map $\psi$.  Notice that the edges of $B_c$ (Figure 15) are exactly the edges which cross the $\tau$-midline in Figure 16.  

\begin{figure}[htb]
\noindent
\includegraphics[width=1.0\textwidth]{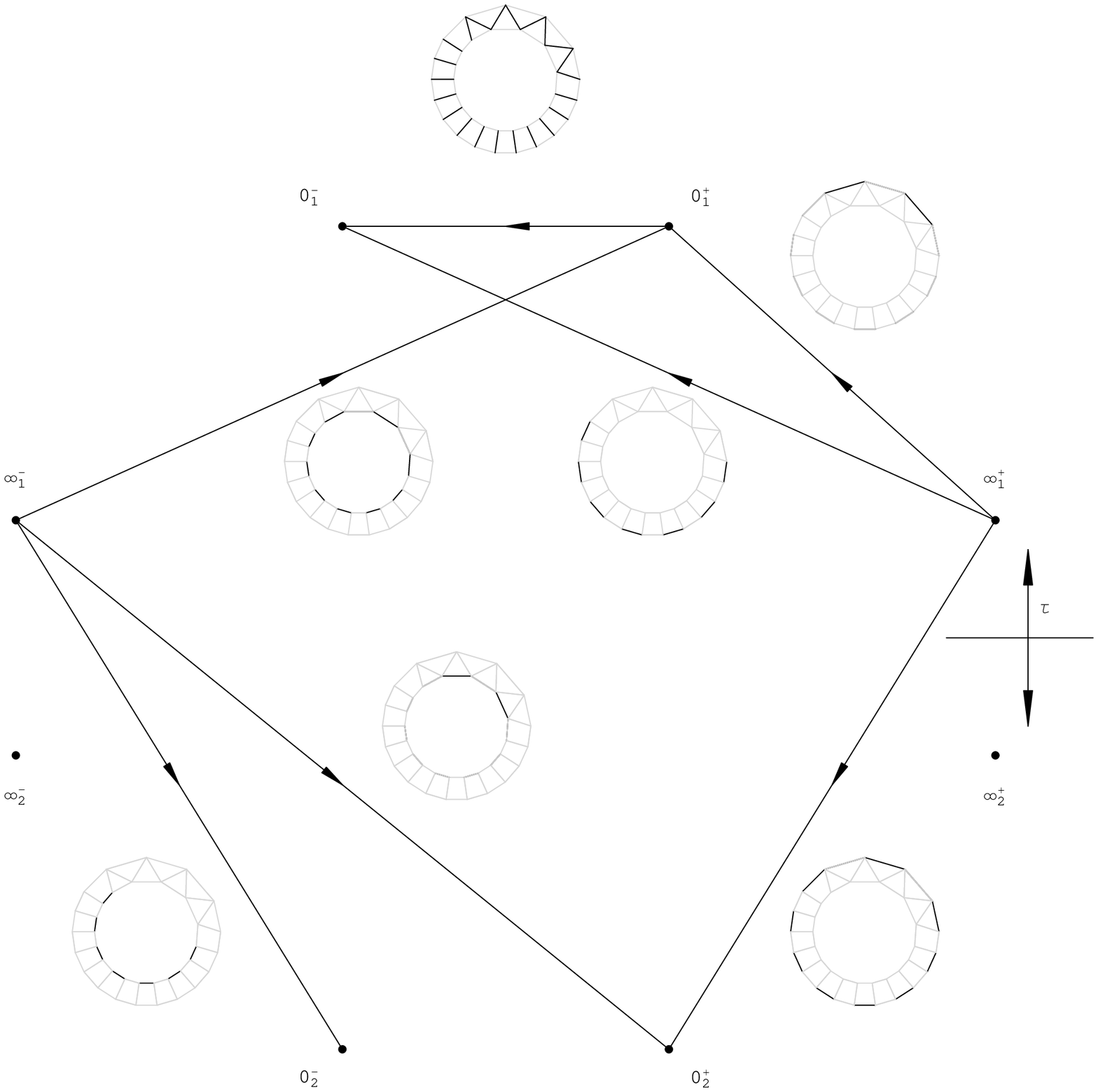}

\ \ \ \ \ \ \ \ \ Figure 16. $\bar\Gamma$.
\end{figure}

We graph the values of $\Psi$ for the four faces in Figures 17-20.

\begin{figure}[htb]
\noindent
\includegraphics[width=.23\textwidth]{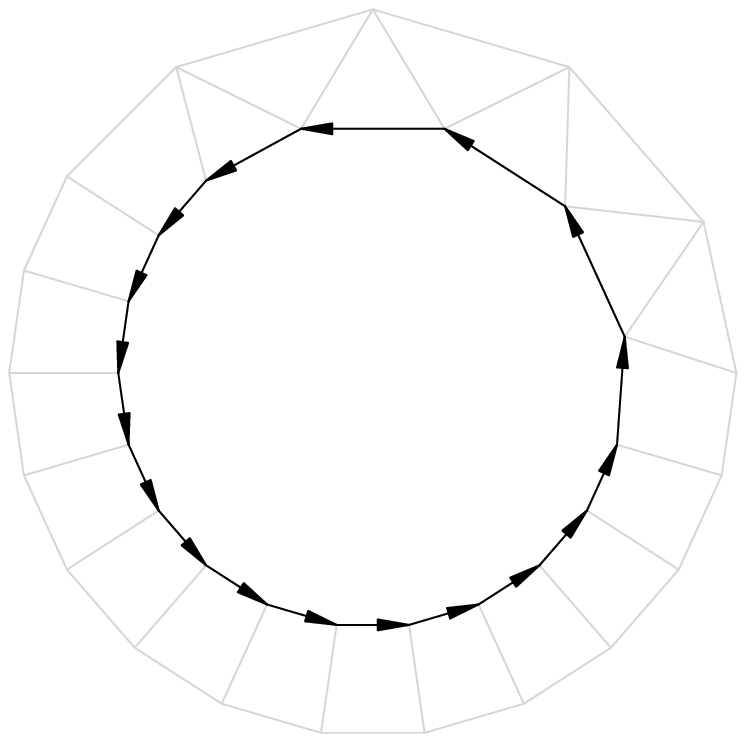}
\ \
\includegraphics[width=.23\textwidth]{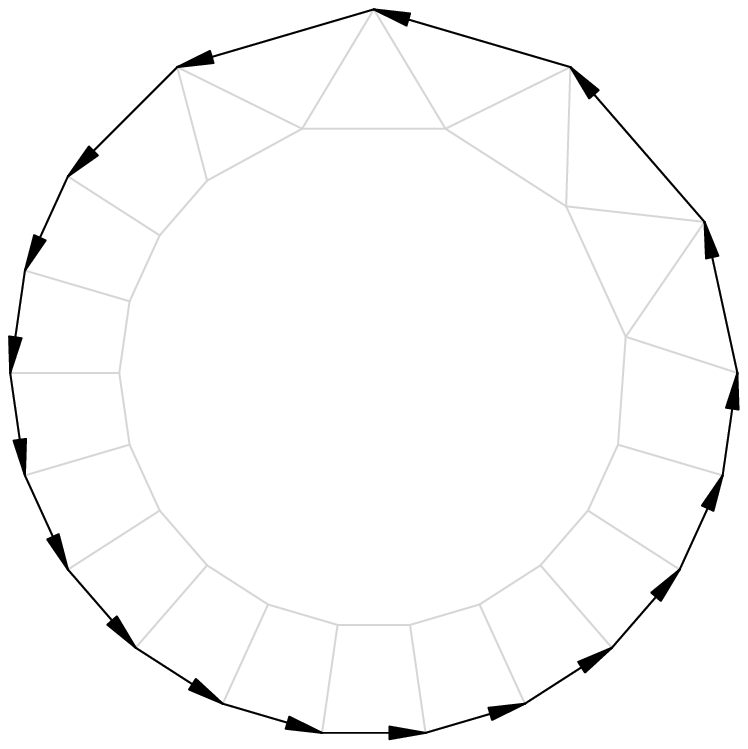}
\ \ 
\includegraphics[width=.23\textwidth]{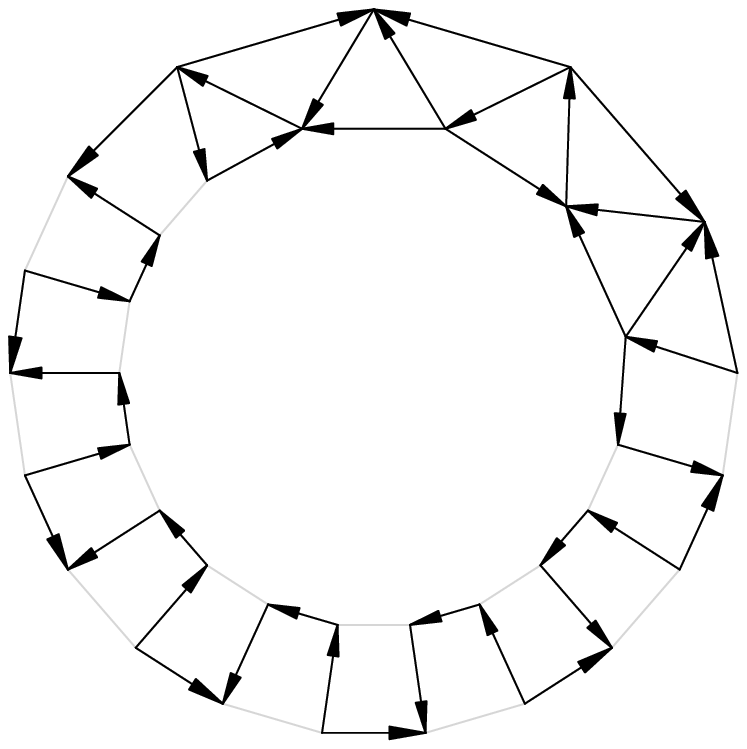}
\ \ 
\includegraphics[width=.23\textwidth]{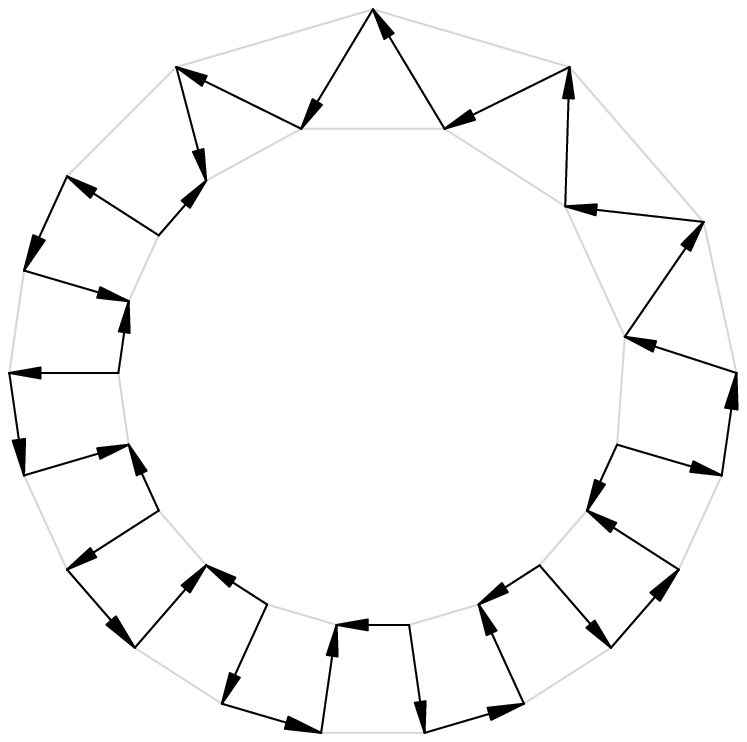}
\noindent
Figure 17. $-\Psi(\infty^-)$\ \ Figure 18. $-\Psi(\infty^+)$\ \ Figure 19. $-\Psi(0^+)$\ \ Figure 20. $\Psi(0^-)$
\end{figure}

Algebraically, we may label our edges counterclockwise from the real axis in sets upper branch ($u_j$), lower branch ($l_j$) or branch locus ($b_j$), depending on where their realizations in $\Gamma$ lie.  We illustrate this (again in genus 10) in Figure 21, from which the pattern can be deduced.   From this we explicitly write $\Psi$ on the faces.  Recall that $\Psi(\mathbb ZF)\subset \psi(\mathbb ZE)$, and $\psi$ is injective, so we will also write $\Psi$ for $\psi^{-1}\circ\Psi:\mathbb ZF\to\mathbb ZE$.

\begin{proposition}\label{prop.ZF}
$\mathbb Z E$ is spanned by elements of the form
\begin{itemize}
\item $u_j$, $1\leq j\leq 2g-2$
\item $l_j$, $1\leq j\leq 2g-2$
\item $b_j$, $1\leq j\leq 2g+2$
\end{itemize}
and we may choose $\psi$ so that a basis for $\mathbb Z F\hookrightarrow^\psi \mathbb Z E$ is given by:
\begin{align}
-\Psi(\infty^-)&=\sum_{j=1}^{2g-2}l_j\notag\\
-\Psi(\infty^+)&=\sum_{j=1}^{2g-2}u_j\notag\\
-\Psi(0^+)&=l_1+l_3-u_2-u_4+\sum_{j=1}^{g+1}(u_{2j-1}-l_{2j})+\sum_{k=1}^{2g+2}b_k\notag\\
\Psi(0^-)&=\sum_{j=3}^{g+1}(u_{2j}-l_{2j-1})+\sum_{k=1}^{2g+2}b_k.\notag
\end{align}
\end{proposition}

These will represent our standard choice for $(\psi,\Psi)$.  We collect this data notationally in the following lemma:

\begin{lemma}\label{lem.boundaries}
The boundaries of the faces of $\cbG$ are cycles joined head-to-tail as follows.  Brackets have been used to signify ``important'' subwords, and carriage returns are placed between the square edges and triangle edges in $\partial 0^{\pm}_1$, as these have different word patterns.  It is understood that these are cyclic words.  For simplicity, for each edge, $e\in E$, we choose some lift $\bar e$ such that $\psi(e)=\bar e-\tau\bar e$.  However, we supress the bar in our notation.

\begin{align}
-\partial \infty^-_1&= l_1l_2l_3\cdots l_{2g-2},\notag\\
-\partial \infty^+_1&= u_1u_2u_3\cdots u_{2g-2},\notag\\
\partial 0^-_1&=[b_1\cdots b_8]\notag\\
&\ \ \ \ \ \ \ \ \ [b_9(\tau l_5)b_{10}u_6]\cdots [b_{2g+1}(\tau l_{2g-3})b_{2g+2}u_{2g-2}],\notag\\
\partial 0^+_1&= [(-b_1)(\tau u_1)(-b_2)(\tau l_1)][(-b_3)u_2(-b_4)l_2][(-b_5)(\tau u_3)(-b_6)(\tau l_3)][(-b_7)u_4(-b_8)l_4]\notag\\
&\ \ \ \ \ \ \ \ \ [(-b_9)(\tau u_5)(-b_{10})l_6]\cdots[(-b_{2g+1})(\tau u_{2g-3})(-b_{2g+2})l_{2g-2}]\notag
\end{align}

\end{lemma} 

\begin{figure}
\noindent
\includegraphics[width=.45\textwidth]{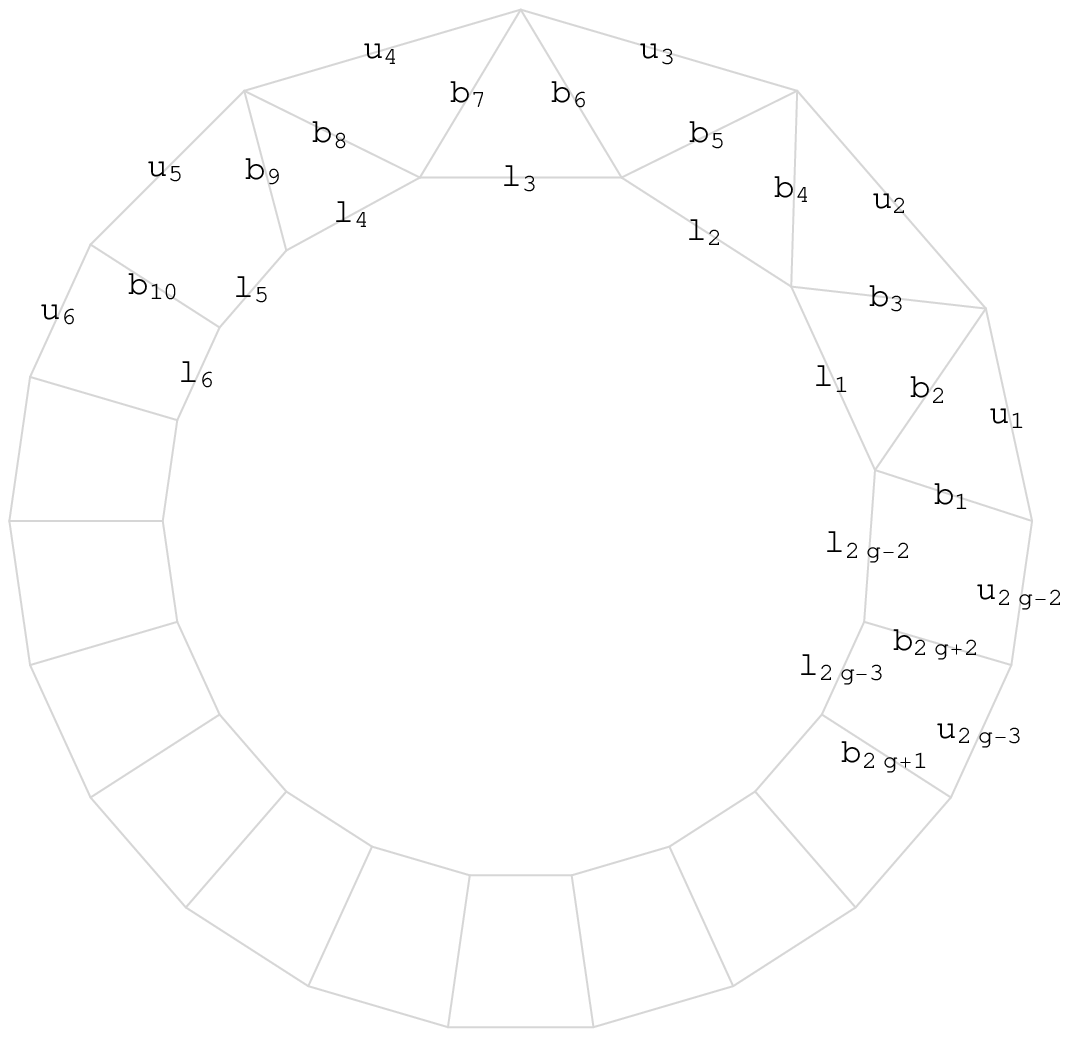}

Figure 21. Enumeration of the edges.
\end{figure}

\section{intersection pairings}
Theorem \ref{thm.repwithint} requires the computation of the intersection pairings on $S$ of the cycles $\psi(e)$.  For completeness, we include these computations in this section.

\begin{lemma}\label{lem.faceedges}
Let $f\in F$ be a face of $\cbG$, and let $e'$, $e$ be oriented edges of $\cbG$ such that $\partial e=x-y$ and $\partial e'=y-z$.  Assume also that $(e',y,e)$ is an interval on $\partial f$.  Then $e-\tau e$ and $e'-\tau e'$, viewed as loops in $S$ have intersection pairing
\[
(e-\tau e)\cdot(e'-\tau e')=+1.
\]
\end{lemma}

\begin{proof}
Notice that $x\neq y$ and $z\neq y$ by construction, and that $e$ and $e'$ only intersect at $y\in S$.  If we look at a neighborhood of $y$ in $S$, We arrive at Figure 22. 

\begin{figure}
\noindent
\includegraphics[width=.3\textwidth]{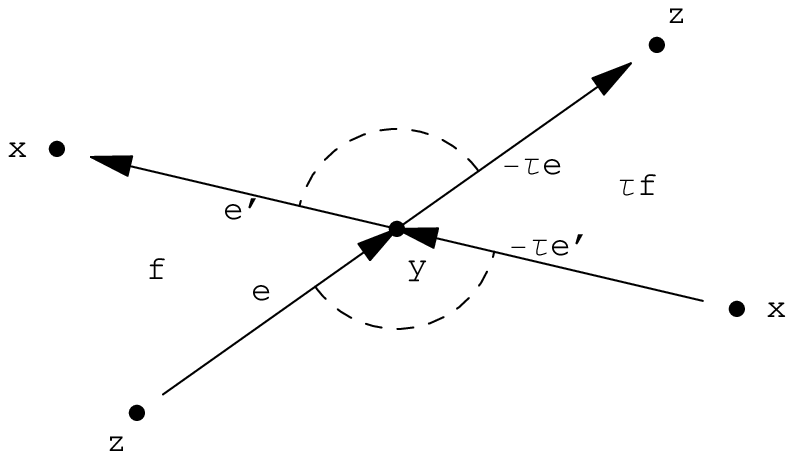}

Figure 22. Computing $(e-\tau e)\cdot(e'-\tau e')$.
\end{figure}

\noindent 
The dotted lines indicate other faces/edges.  From this figure, the result is clear.
\end{proof}

Thus, for example, if we let $f=\infty^+$, $e=-u_2$, $e'=-u_1$, so that 
\[
(-u_2+\tau u_2)\cdot (-u_1+\tau u_1)=+1.
\]
This lemma in conjunction with Lemma \ref{lem.boundaries} gives the intersection numbers:

\begin{proposition}\label{prop.intnos}
Let $u_j$, $l_j$, and $b_j$ denote the cycles in $S$ labelled as in Figure 21.  Notice that we are now using the symbols $e$ to denote $\psi(e)=e-\tau e$.  Let $e$, $e'$ denote any two of these cycles.
\begin{enumerate}
\item $e\cdot e'=0$ if $e$ and $e'$ do not share a vertex.
\item $u_j\cdot u_{j+1}=-1$
\item $l_j\cdot l_{j+1}=-1$
\item $b_j\cdot b_{j+1}=+1$ for $1\leq j\leq 7$
\item $u_j\cdot b_{2j}=-u_j\cdot b_{2j+1}-u_{j+1}\cdot b_{2j}=u_{j+1}\cdot b_{2j+1}=+1$, for $j=1,3$
\item $u_j\cdot b_{2j}=-u_j\cdot b_{2j+1}-u_{j+1}\cdot b_{2j}=u_{j+1}\cdot b_{2j+1}=-1$, for $j=2,4$
\item $l_{j-1}\cdot b_{2j-1}=-l_{j-1}\cdot b_{2j}=-l_{j}\cdot b_{2j-1}=l_{j}\cdot b_{2j}=-1$, for $j=1,3$
\item $l_{j-1}\cdot b_{2j-1}=-l_{j-1}\cdot b_{2j}=-l_{j}\cdot b_{2j-1}=l_{j}\cdot b_{2j}=+1$, for $j=2,4$
\item $l_j=b_{j+5}=-l_{j+1}=b_{j+5}=-1$, for $4\leq j\leq 2g-3$
\item $u_j=b_{j+5}=-u_{j+1}=b_{j+5}=+1$, for $5\leq j\leq 2g-2$
\end{enumerate}
\end{proposition}

\begin{proof}
(1) follows from noticing that $\check\Gamma$ may be realized inside $M$, so if $e$ and $e'$ do not intersect on $\check\Gamma$, then they lift to two disjoint sets on $S$.  

If $e$ and $e'$ are consecutive edges of a word in Lemma \ref{lem.boundaries}, then $(e-\tau e)\cdot (e'-\tau e')=+1$ by Lemma \ref{lem.faceedges}.  This proves the result for all cases except if $e$ and $e'$ are opposite edges on a four-valent vertex.  However, if $ee''$ is a consecutive pair of edges on $f_a$ and $(e'')^{-1}e'$ is a pair of edges on $f_b$, then we may apply the argument of Lemma \ref{lem.faceedges} to $f_a\cup_{e''} f_b$, which has consecutive edges $ee'$, to arrive at the result.
\end{proof}

\section{The classical Burau representation}
In this section, we recall the classical Burau representation, and discuss a specialization, which may be helpful later as an analogy.  One may find a detailed description in \cite{Birman}.

Throughout this section, if $R$ is any Riemann surface, then let $\sigma\mapsto\bar\sigma$ be the Hurewicz homomorphism, $\pi_1(R)\to H_1(R,\mathbb{Z})$.  Also use the same notation for the map $\pi_1(R^{[n]})\to H_1(R,\mathbb{Z})$ defined implicitly\footnote{That is, defined only for the case $n=4g-4$.} in Proposition \ref{prop.inj}.

We construct the classical Burau representation.  Let $B_n=\pi_1(D^{[n]})$ be the braid group on the disk, $D$, $X\in D^{[n]}$.  Let $w$ be the winding number function $w:H_1(D\setminus X,\mathbb Z)\to \mathbb Z$ such that $w(\partial x)=+1$ for all $x\in X$.  $w$ induces a covering space $\tilde D\to D\setminus X$ by the map
\begin{align}
\pi_1(D\setminus X)&\to \mathbb Z\notag\\
\sigma&\mapsto w(\bar\sigma).\notag
\end{align}
This space is equipped with a vertical translation operator $t$.  This is the map induced by moving the basepoint ``up'' along a curve of winding number one, thus increasing the winding number of all points in $\tilde D$ by one.  Let $\Lambda=\mathbb Z[t,t^{-1}]$.  $H_1(\tilde D,\mathbb Z)$ is a finite rank $\Lambda$-module, and gives a $\Lambda$-representation $B_n\to \text{Aut}_\Lambda(H_1(\tilde D,\mathbb Z))$.  This is the Burau representation.  

Explicitly, $\pi_1(D^{[n]})$ is generated by transpositions $\sigma_1,\ldots,\sigma_{n-1}$, and the Burau representation is an action on $\Lambda X$, given by 
\begin{align}
\sigma_jx_j&=(1-t)x_j+tx_{j+1}\notag\\
\sigma_jx_{j+1}&=x_j\notag\\
\sigma_jx_k&=x_k\text{~~~~otherwise}\notag
\end{align}
One finds these relations by letting $x_j$ represent a curve from $\partial D$ to $\partial D$ with only $x_j$ in its interior, lifting to $\tilde D$, and rescaling appropriately by some power of $t$.  

We consider two specializations, which we will generalize to surfaces.  Consider the quotient $t^k=1$.  This is equivalent to constructing the $k$-fold cover $\tilde D_k\to D\setminus X$ given by 
\[
w~~\text{mod~} k:\pi_1(D\setminus X)\to \mathbb Z/k\mathbb Z.
\]
Observe that if $k=2$, and we let $\xi_j=x_j+x_{j+1}$, then $\xi_j$ is homotopy equivalent to a loop in $\tilde D_2$ around the points $x_j,x_{j+1}\in X$.  Morally these are the cycles $e-\tau e$.  On $\tilde D_2$, the nontrivial intersection pairings are 
\[
\xi_j\cdot \xi_{j+1}=-1.
\]
One can compute directly from the definitions above that 
\[
\sigma_j(\xi_k)=\xi_k-(\xi_k\cdot \xi_j)\xi_j.
\]

Compare this result with Theorems \ref{thm.teaser} and \ref{thm.repwithint}.  Indeed, one can apply exactly the argument from Theorem \ref{thm.repwithint} to prove this as well.

Now consider the quotient of the Burau representation by $1+t+t^2+\ldots+t^{k-1}=0$.  To realize this we construct the cover $\bar D_k\to D$ which is the $n$-point compactification of $\tilde D_k$ by the ramification points, $X$.  This is distinguished from $\tilde D_k$ by the fact that the loop around the point $x_j\in X$, which is $(1+t+t^2+\ldots+t^{k-1})x_j$ is now contractible to 0.  We call this the compact Burau representation associated to the quotient $\Lambda/t^k-1$, to be consistent with the following section.  This has no counterpart for the Burau representation, since 
\[
\Sigma_{i\in \mathbb Z}t^i\notin\Lambda.
\]

We will mainly be interested in the case $k=2$, so the quotients $t^2=1$ and $t=-1$.

Observe that we can construct these quotients of the Burau representations as monodromy representations of bundles.  Let 
\[
\bar{\mathcal D}_k=\{(X,x,y)\in D^{[n]}\times D\times \mathbb C|y^k=\prod_{x_j\in X}(x-x_j)\},
\]
and let $\mathcal D_k=\{(X,x,y)\in\bar{\mathcal D}|y\neq 0\}$.  Then the fibers of these bundles are $\bar D_k$ and $\tilde D_k$, and the monodromy representations are the quotients discussed above.

\section{The surface Burau representation}\label{sec.burau}

In this section, we construct a generalization of the classical Burau representation.  Later we show that, up to twisting by $\tau$, the homology of $\mathcal M$ is related to the compact form of this Burau representation via pullback from $M^{[4g-4]}$.  This will tell us the monodromy representation of $\pi_1(H^s)$.  

There are three obstructions to generalizing the Burau representation to a general Riemann surface $M$ of positive genus: $w$ is well-defined only up to $\mathbb Z/n\mathbb Z$, there is no \textit{a priori} choice for $w$ on the complement of $\text{ker}(H_1(M\setminus X,\mathbb Z)\to H_1(M,\mathbb Z))$, and monodromy would be determined only up to automorphisms of the cover (as there is no boundary to use to fix the identification). 

Let $M^{[n]+1}=\{(m,X)\in M\times M^{[n]}|m\notin X\}$.  This is consistent with the notation of Section \ref{sec.surf}.  For $(m,X)\in M^{[n]+1}$, some choice of winding number, $w$, and some positive integer $k| n$, we will construct a $k$-fold covering of $M\setminus X$ with a marked point.   Let $(B',m')$ be the based covering space of $(M\setminus X,m)$ associated to the map $\pi_1(M\setminus X)\to \mathbb Z/k\mathbb Z$ given by ${\sigma\mapsto w(\bar\sigma)\mod k}$.  We need the restriction $k|n$ since the ``interior'' of $\bar\sigma$ is not a well-defined notion.\footnote{$w_{\text{int}}(\bar\sigma)-w_{\text{ext}}(\bar\sigma)=\pm n$.}
We may complete this space by adding $X$ as ramification points.  Call this new space $B^k_{mXw}$.  It is topologically a surface of genus 
\[
\tilde{g}=k(g-1)+1+\frac{(k-1)n}{2}.
\]
This space has a vertical translation operator, $t$, so its homology is a 
\[
\Lambda_k=\mathbb{Z}[t]/\langle t^k-1\rangle\]
module.  Let $\tilde B^k_{mXw}=B^k_{mXw}\setminus X$, so that $B^k_{mXw}$ is the compact form (compactification by $X$) of $\tilde B^k_{mXw}$.

\begin{definition}
$\WW\to M^{[n]+1}$:
\end{definition}

Let $W_X^k$ be the set of winding number functions on $M\setminus X$ modulo $k$:
\[
W_X^k=\{w\in H_1(M\setminus X,\mathbb Z/k\mathbb Z)|w(\partial x)=+1\forall x\in X\}.
\]
This is a set with $k^{2g}$ points.  Notice that this is affine, modelled on the lattice $H_1(M,\mathbb Z/k\mathbb Z)^*$.  Let $X\in M^{[n]}$, and choose $U\subset M$ contractible open such that $X\subset U$.  Then $U^{[n]}$ gives a neighborhood of $X$ in $M^{[n]}$.  A winding number function on $M\setminus X$ is defined entirely by its restriction to $H_1(M\setminus U, \mathbb Z)$.  Thus $w$ is well-defined as an element of $W_Y^k$ for all $Y\in U^{[n]}$, giving natural isomorphisms $W_Y^k\cong W_X^k$.  We have bundles $\bar\WW|_{U^{[n]}}\to U^{[n]}$ for all contractible open $U\subset M$ with fibers $W_Y^k$, and restriction to $U'\subset U$ gives a canonical isomorphism 
\[
(\bar\WW|_U)|_{U'}\simeq \bar\WW|_{U'},
\]
so this gives a well-defined locally constant bundle $\bar\WW\to M^{[n]}$.
Let 
\[
\WW=\bar\WW\times_{M^{[n]}}M^{[n]+1}.
\]

\begin{definition}
$\BB\to\WW$:
\end{definition}

Let $\BB$ be the set 
\[
\BB=\coprod B^k_{mXw}.
\]

$\BB$ has a unique and universal structure as a bundle over $\WW$:

\begin{lemma}\label{lem.universal}
$\BB$ satisfies the following properties:

\begin{itemize}
\item $\BB$ is a topological bundle over $\WW$.
\item If $\WW'\subset \WW$ and $\BB'\to\WW'$ is another bundle such that $\BB'|_{(m,X,w)}=B^k_{mXw}$ with section $\mu$ such that $\mu(m,X,w)$ lies over $m$ (a section of marked points), then there exists a unique isomorphism $\BB'\simeq \BB|_{\WW'}$ preserving the marked point sections.
\end{itemize}
\end{lemma}

\begin{proof}
Fix $(m,X,w)\in \WW$.  Let $U,V\subset M$ contractible, $U\cap V=\emptyset$, $X\subset U$, $m\in V$.  This defines an open neighborhood of $(m,X,w)$ in $\WW$, since $\WW$ is discrete over $M^{[n]+1}$.  The topology of $\WW$ is generated by such neighborhoods.  For any $m\in V$, $w$ gives a well-defined unique marked double cover of $M\setminus U$, and all are canonically isomorphic, as $V$ is contractible.  Now one uses the classical Burau construction to build the double cover on $U$, and identify the boundary.  Specifically, choose a smooth coordinate $z:U\to \mathbb C$ and construct the cover  
\[
\{(y,x)\in\mathbb C\times U|y^2=\prod_{x_j\in X} (z-x_j)\}\to U^{[n]}\times V.
\]
Since the order of $X$ is even, the boundary of this cover is a disjoint union of circles on which monodromy of $\pi_1(U^{[n]})$ acts trivially, so the gluing construction is well-defined.

We now need to show that under restriction in the open cover:
\[
\left(\BB|_{(V,U^{[n]},w)}\right)|_{(V',(U')^{[n]},w)}=\BB|_{(V',(U')^{[n]},w)}
\]
for $V'\subset V$ and $U'\subset U$, but this must be true, as on the overlap $U\setminus U'$, each is just the disjoint union of two annuli: $U\setminus U'\times \{0,1\}\times(U')^{[n]}\times V'$, and elsewhere canonically isomorphic, given that we choose the same coordinate on $U$ and $U'$.  The first point is proved.

Now assume that there is some $\BB'\to \WW'$, with marked point $\mu$.  There is a canonical map to $\BB|_{\WW'}$ by the uniqueness of $B^n_{mXw}$.  We need to show this map is continuous.  However, one can recover from $\WW'$ and $\mu$ the geometry of $\BB$ via the exact construction that gave the first part of the lemma.  Then the bundles are easily locally isomorphic (via the global map), so are isomorphic.
\end{proof}

We set out to construct a bundle over $M^{[n]}$ whose monodromy generalizes the Burau representation.  We have arrived at: $\BB=\coprod B^n_{mXw}\to\WW\to M^{[n]+1}\to M^{[n]}$.  This is a far cry from what we want, as it has a legion of fibers telling the extra data of winding numbers and marked points.  Now we try to reduce this information.  First we study the monodromy of $\WW\to M^{[n]+1}$.

\begin{proposition}
Let $\langle\cdot,\cdot\rangle$ denote the intersection pairing on $H_1(M)$.   
\begin{itemize}
\item Every element of $\pi_1(M^{[n]+1},(X,m))$ may be factored as $\sigma\sigma'$, for \\$\sigma\in \pi_1((M\setminus m)^{[n]},X)$ and $\sigma'\in\pi_1(M\setminus X,m)$
\item $\sigma'(w)=w$ for  $\sigma'\in \pi_1(M\setminus X)$
\item $\sigma(w)=w-\langle \bar\sigma,\cdot\rangle$ for $\sigma\in \pi_1((M\setminus m)^{[n]})$
\end{itemize}
\end{proposition}

\begin{proof}

Applying Corollary \ref{cor.fn2}, we see that 
\[
\pi_1((M\setminus m)^{[n]})\times\pi_1(M\setminus X)\twoheadrightarrow\pi_1(M^{[n]+1}).
\]

If $\sigma'\in\pi_1(M\setminus X)$, then $\sigma'$ acts trivially on $\WW$, since $\WW=\WW'\times_{M^{[n]}}M^{[n]+1}$.  

Assume  $\sigma\in \pi_1((M\setminus m)^{[n]})$.  
Let $a_j$, $b_j$ be a symplectic basis of $H_1(M,\mathbb Z)$ based at $x_1$ and missing $m$.   
Since $\pi_1(M\setminus X)$ acts trivially and $\sigma\to\bar\sigma$ is a group homomorphism, we may assume by Theorem \ref{thm.structure} that $\sigma$ is one of the following:
\begin{itemize}
\item a transposition of $x_j$ and $x_k$ across some edge $e$
\item $\sigma$ fixes $x_j$, $j>1$, and $\sigma$ follows $a_j$ or $b_j$.  
\end{itemize}
We study $\sigma(w)-w$:
\[
\sigma(w)(\alpha)-w(\alpha)=w(\sigma^{-1}\alpha-\alpha).
\]

Assume that $\sigma$ is a transposition around an ``edge,'' $e$ connecting $x_1$ to $x_k$.  Then $\sigma^{-1}\alpha-\alpha$ is supported in a neighborhood of $e$, so
\[
\sigma^{-1}\alpha-\alpha=a\partial x-b\partial y.
\]
However, applying $\sigma$ (an involution) to both sides, we see
\[
\alpha-\sigma^{-1}\alpha=-b\partial x+a\partial y,
\]
so $a=b$, and 
\[
w(\sigma^{-1}\alpha-\alpha)=w(a\partial x-a\partial y)=a-a=0.
\]
The propositon follows in this case, since $\bar\sigma=e-e=0$.

Now assume that $\sigma$ fixes $x_j$, for $j>1$, and $x_1$ traces $a_k$.  Again we study $\sigma^{-1}\alpha-\alpha$.  This is zero for $\alpha=b_j$ for $j\neq k$, $\alpha=\partial x_j$, $j>1$, and  $\alpha=a_j$, since the supports are disjoint ($\sigma$ may be realized as an automorphism fixing all but a neighborhood of $a_j$).  Also, 
\[
\sigma\partial x_1=\partial\sigma x_1=\partial x_1,
\]
So $\sigma$ is trivial on all basis elements but $b_k$.  

Figure 23 illustrates how $\sigma b_j$ acquires $\partial x_1$ as $x_1$ passes through $b_j$.  Replacing $(a_j,b_j)$ by $(b_j,-a_j)$ shows the result for $\bar\sigma=b_j$ as well.  Thus 
\[
\sigma w(\alpha)=w(\alpha)-\langle\bar\sigma,\alpha\rangle
\]
for all generators $\sigma$ of $\pi_1((M\setminus m)^{[n]})$.  

\begin{figure}[htb]
\noindent
\includegraphics[width=.3\textwidth]{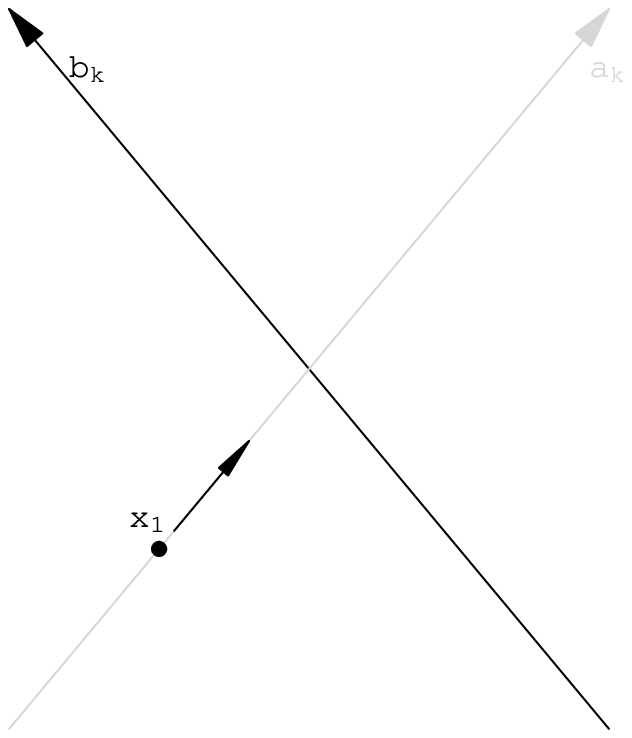}~~~~~~~~
\includegraphics[width=.3\textwidth]{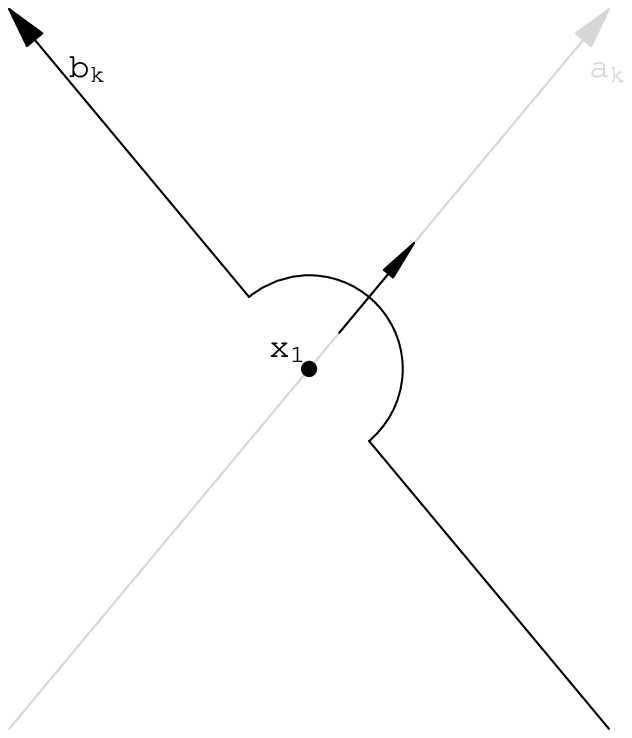}~~~~~~~~
\includegraphics[width=.3\textwidth]{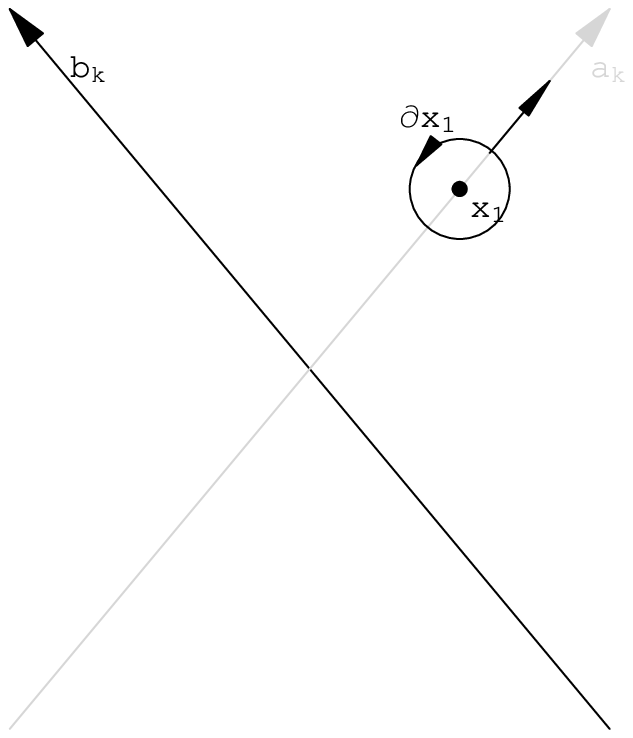}

Figure 23.  $\sigma b_k$.
\end{figure}
\end{proof}

If $\sigma\in\pi_1(M\setminus X)$, then $\sigma:\BB|_{(m,X,w)}\to\BB|_{(m,X,w)}$ acts by vertical translation by $t^{w(\bar\sigma)}$ (specifically, it preserves the components of $\BB\to M^{[n]+1}$) since $\sigma$ moves the basepoint.  If $\sigma\in\pi_1((M\setminus m)^{[n]})$, then $\sigma$ acts trivially on $\WW$ (preserves the components) only when $\bar\sigma\in kH_1(M\setminus X,\mathbb Z)$.  This is because $\sigma(w)=w+\langle\cdot,\bar\sigma\rangle$.

For now we fix $k$ and let $K=\{\sigma\in\pi_1((M\setminus m)^{[n]})|\bar\sigma\in \mathbb{Z}\partial m\}$.  The group $K$ is the kenel of 
\[
\pi_1((M\setminus m)^{[n]})\to H_1(M,\mathbb Z).
\]
Then for any $w\in W_X$, $\sigma(w)=w$, and $\BB$ gives a monodromy action of $K$ on $B=B^k_{mXw}$.  $H_1(B,\mathbb{Z})$ is a $\Lambda_k$-module, and a $\Lambda_k$-representation of $K$.  Furthermore, if $\sigma\in\pi_1(M\setminus X)$, then $\sigma$ acts by $t^{w(\bar\sigma)}$ as it changes only the basepoint.  Thus in fact, we get a map:
\[
\text{ker}(\pi_1(M^{[n]})\to H_1(M,\mathbb Z))~~\to~~Aut_{\Lambda_k}H_1(B,\mathbb{Z}))/\langle t\cdot \text{Id}\rangle,
\]
which we will call the Burau representation.  This is not a linear representation, and it depends on the choice of $w$.  

Notice that this is $\text{Aut}_{\Lambda_k}$ and not just $\text{Aut}_{\mathbb Z}$, since by Corollary \ref{cor.fn1}, $\langle t\cdot \text{Id}\rangle$ is normal in the image of $\pi_1(M^{[n]+1})\to \text{Aut}_{\mathbb Z}(H_1(B,\mathbb Z))$.

\section{The bundles $\SS$ and $\HH_1(\text{Prym}~\SS)$}\label{sec.bundleS}
At this point, we begin to study structure of the Hitchin bundle.  Letting $\HH_1(\MM)$ denote the bundle of integral first homology of the fibers of $\MM\to H^s$ we are studying the action of $\pi_1(H^s)$ on $\HH_1(\MM)$.  For a construction of $\HH(\cdot)$, see \cite{Griffiths}.  First, we construct the bundle $\SS$.  
We use this bundle to model $\HH_1(\mathcal M)$ by Theorem \ref{thm.MtoJ}.  In Proposition \ref{prop.withdiagram}, we show that $\HH_1(\SS)$ is closely akin to the Burau bundle and prove that $\text{ker}(\pi_1(H^s)\to\pi_1(M^{[4g-4]}))$ acts almost trivially in monodromy (in fact by the group $\{1,\tau\}$).  We finish by studying the monodromy via the combinatorial techniques of Section \ref{sec.phi}.

Henceforth, let $T$ be the total space map, so that $TK_M$ is the total space of the canonical bundle, and $TK^2_M$ the total space of the square of the canonical bundle.
The map $H^s\times M\to TK^2_M$, $(\omega,m)\mapsto\omega(m)$ is holomorphic.  Consider the (holomorphic) map, $H^s\times TK_M\to TK^2_M$ defined by $(\omega,m,q)\to (m,\omega(m)-q^2)$.  We have a variety:

\begin{definition}
The space, $\SS$, is defined as
\[
\SS=\{(\omega,m,q)\in H^s\times TK_M|\omega(m)-q^2=0\}.
\]
\end{definition}

Notice that for $\omega\in H^s$, $\SS|_{\{\omega\}}\cong S_\omega$, naturally.  For any $\omega\in H^s$, there is some neighborhood $U\subset H^s$ of $\omega$ such that $\SS|_U\cong U\times S_\omega$ as smooth manifolds.  Thus $\SS$ is a smooth bundle $\SS\to H^s$ fibered by double covers of $M$.

\begin{definition}
Let $\HH_1(\SS,\mathbb{R})\to H^s$ be the real bundle defined fiberwise by:
\[
\HH_1(\SS,\mathbb R)|_{\{\omega\}}=\HH(S)|_\omega\otimes_{\mathbb Z}\mathbb R.
\]
Continuous local sections descend from $\mathcal{O}_{\mathbb{R}}\otimes \HH(S)$.
Let $\text{Jac}(\SS)\to H^s$ be the bundle defined fiberwise by:
\[
\text{Jac}(\SS)=\HH_1(\SS,\mathbb R)/\HH(S).
\]
\end{definition}

\section{Some differential geometry}\label{sec.var}

Throughout this section, we use $H^1_{DR}$ to denote complexified deRham cohomology.  We introduce a local coordinate system on $\text{Jac}(\SS)$ via integrals in order to return to a classical definition of the Jacobian as a space of integrals.

\begin{proposition}
If $dx_1,\ldots,dx_{\tilde g}\in H^1_{DR}(S_\omega)$ with $dx_1,\ldots,dx_{\tilde g},d\bar x_1,\ldots,d\bar x_{\tilde g} $ linearly independent, then there exist contractible open $U$ containing $\omega$, and $dz_1,\ldots,dz_{\tilde g}\in H^1_{DR}(\SS|_U)$, such that $dz_j|_{S_\omega}=dx_j$ and $dz_j|{\eta}, \overline{dz_j|_\eta}$ linearly independent.
\end{proposition}

\begin{proof}
Choose any such $U$ contractible and open.  Then $ H^1_{DR}(\SS|_U)\cong H^1_{DR}(S_\omega)$ by restriction, so there exist $dz_j$ such that  $dz_j|_{S_\omega}=dx_j$.  Let $\omega$ and $\eta$ be elements of $U$.  Let $\alpha_\omega\in H_1(S_\omega,\mathbb Z)$.  then $\alpha$ is homotopic to a curve $\alpha_\eta\in H_1(S_\eta,\mathbb Z)$, and in fact, $U$ gives a natural isomorphism between the different groups $H_1(S_\eta,\mathbb Z)$ for $\eta\in U$.  Since the $dx_j$ are closed, 
\[
\int_{\alpha_\omega}dx_j=\int_{\alpha_\eta}dx_j,
\]
so that the $dx_j$ give an isomorphism between groups $H^1_{DR}(S_\omega)$ and  $H^1_{DR}(S_\eta)$.  Specifically, the restrictions of the one-forms remain linearly independent.
\end{proof}

In fact we have the much stronger corollary, which is a family version of Abel-Jacobi (notice that this is not a holomorphic statement):

\begin{theorem}\label{thm.trivialjac}
If $U$ is contractible and open, then there exists $dx=(dx_1,\ldots,dx_{\tilde g})$, $dx_j\in H^1_{DR}(\SS|_U)$ such that $\int_\bullet dx$ gives an isomorphism:
\[
\text{Jac}(\SS)|_U\to U\times \mathbb C^{\tilde g}/\Lambda,
\]
where $\Lambda=\{\int_\alpha dx|\alpha\in H_1(S_\omega,\mathbb Z)\}$ for some (any) $\omega\in U$.  
\end{theorem}

\begin{proof}
We need only show that if $\alpha\neq 0\in H_1(S_\eta)\otimes \mathbb R$, then $\int_\alpha dx\neq 0$.  However, if $\int_\alpha dx= 0$, then $\int_\alpha d\bar x= 0$, and the $dx_j$, $d\bar x_j$ form a basis for complexified deRham cohomology, which gives a contradiction.
\end{proof}

Given this, we can treat degree 0 divisors on $S_\omega$ as elements of $\text{Jac}(\SS)|_{\{\omega\}}$.  If $D\in \text{Div}^0(S_\omega)$, then 
\[
\int_D dx\equiv\int_{\alpha_D}dx\mod \Lambda
\]
for some $\alpha_D\in H_1(S_\omega,\mathbb Z)\otimes_{\mathbb Z}\mathbb R$.  However, since $dx_j$, $d\bar x_j$ span $H^1_{DR}(S_\omega)$, $D\mapsto \alpha_D$ is independent of coordinate $dx$ and is well-defined as a map from $\text{Div}^0(S_\omega)$ to $\text{Jac}(\SS)|_{\{\omega\}}$.  

\begin{theorem}
If $U$ is contractible, $dx$ is as above, and $d_1,\ldots,d_n, e_1,\ldots,e_n$ are smooth sections of $\SS|_U$, then 
\[
\sum_{j=1}^n d_j(\omega)-e_j(\omega)\in \text{Div}^0(S_\omega)
\]
for all $\omega\in U$, and $\sum_j\int_{e_j(\omega)}^{d_j(\omega)}$ is a smooth section of $\text{Jac}(\SS)|_U$.
\end{theorem}

\begin{proof}
By Theorem \ref{thm.trivialjac}, we need only show that
\[
\omega\to \sum_j\int_{e_j(\omega)}^{d_j(\omega)}dx
\]
is smooth, but this is clear.
\end{proof}

\section{A canonical section and the isomorphism $\HH_1(\mathcal M)\to\HH_1(\SS)^{-}$.}
In this section, we finally prove that $\HH_1(\MM)$ is modelled on $\HH_1(\SS)$ (Theorem \ref{thm.MtoJ}).  We have chosen this differential geometry approach in order to hide the most difficult mathematics in Lemma  \ref{lem.veryhard}.  This uses a deep result about index theory to show that some line bundle has a ``meromorphic'' section, as one might expect.  The casual reader may accept this fact and pass by.  One might also choose to consider the isomorphism from a holomorphic point of view, but this leads quickly to a discussion of coherent sheaves.

Let $\pi_\omega:S_\omega\to M$.  To each $\omega\in H^s$, we construct a canonical section $\sqrt{\omega}\in H^0(S_\omega,\pi_\omega^*K_M)$.  We have the line bundle $\pi_\omega^*K_M\to S_\omega$.  A section of this map $\gamma:S_\omega\to\pi_\omega^*K_M$ is an assignment $\gamma(s)\in K_M|_{\pi(s)}$.  Let $\gamma(s)=s$ be our section.  More formally, the identity map $S_\omega\hookrightarrow TK_M$ factors through $T(\pi^*_\omega K_M)$, and $S_\omega\to T(\pi^*_\omega K_M)$ is a natural, single-valued section of $\pi_\omega^*K_M$.  Call this section $\sqrt{\omega}$.

$\SS\subset H^s\times TK_M$, so let $\Pi:\SS\to H^s\times M$ be the projection on the second factor.  Note that $\Pi|_{\{\omega\}}=\text{id}\times \pi_\omega:\{\omega\}\times S_\omega\to\{\omega\}\times M$.  

The map $\text{id}\times (-\text{id}):H^s\times TK_M\to H^s\times TK_M$ restricts to the involution $\tau$ on $S_\omega$, so that $\Tau=\text{id}\times (-\text{id})|_{\SS}$ is the involution such that $\Tau|{\{\omega\}}=\tau$.

Recall that 
\[
\mathcal{M}=\{(\Phi,V)\text{~stable pairs,~}|\wedge^2V=\xi
\text{~for some fixed~}\xi\}.  
\]
V is holomorphic, so compatible with some operator $\bar\partial$ defining the holomorphic structure.

\begin{lemma}\label{lem.linebun}
For $\omega\in U\subset H^s$ contractible open, there exists $U'\subset U$ contractible open containing $\omega$ and a line bundle $\LL\to \SS|_{U'}$ of degree 1 such that $\LL|_{S_\omega}$ is holomorphic, and $\LL$ has a section, $s$ such that $s|_{S_\omega}$ is not the zero section for $\eta\in U'$.
\end{lemma}

\begin{proof}
Choose any degree one holomorphic line bundle. 
\end{proof}

In the next lemma, we show that for any family of Riemann surfaces with line bundle, locally there is a family of nontrivial meromorphic sections.  This relies on an index theory result about the family of complex structures $\bar\partial$.

\begin{lemma}\label{lem.veryhard}
Let $U\subset H^s$ be a contractible open set containing some $\omega\in H^s$ and let $\LL\to \SS|_{U'}$ be a line bundle which is holomorphic on the fibers ($\bar\partial_\omega\LL=0$ for all $\omega\in U'$).  
Then there exists some line bundle $\LL'$ and sections $s:\SS|_U\to \LL\otimes \LL'$ and $s':\SS|_U\to\LL'$ with $s|_{S_\eta}\neq 0$ for $\eta$ in some open subset $U'\subset U$ containing $\omega$.
\end{lemma}

Effectively, this says that $s/{s'}$ is a ``meromorphic'' section of $\LL$.  We won't use any properties of $H^s$ in the proof.

\begin{proof}
Let $U$, $\LL_1$, and $s_1$ be given by Lemma \ref{lem.linebun}.  Let $d$ be the degree of $\LL|_{S_\omega}$ (which is independent of $\omega$).  Let $k$ be any positive integer such that $d+k> 2\tilde g-2$.  Let $\LL'=\LL_1^{\otimes k}$ and $s'=s_1^k$.  Then the line bundle $\LL\otimes \LL'$ is a smooth line bundle, holomorphic on each $\LL\otimes \LL'|_{S_\omega}$.  This is a family of holomorphic line bundles, equipped with a family of elliptic operators, $\bar\partial_\omega$.  Since the degrees at each $\omega$ are larger than $2\tilde g-2$, these operators are of constant positive rank.  Thus there is a family of solutions over $U$ (contractibility), which is the desired section, $s$.
\footnote{This nontrivial result may be found as Theorem 9.11 in \cite{Berline}.}
We may assume this family is nonzero on some restriction of $U$.
\end{proof}

\begin{proposition}\label{prop.transportsec}
Let $\LL\to \SS|_U$ be a line bundle holomorphic and degree zero along the $S_\omega$, $U$ contractible.  Then there exists a section $D:U\to\text{Jac}(\SS)|_U$ such that $\LL|_{S_\omega}=[D(\omega)]$.
\end{proposition}

\begin{proof}
Applying the lemma, we get $s$, $s'$ such that $s/s'$ is a ``meromorphic'' section of $\LL$ over some subset $U'\subset U$.  Then the family of divisors associated to $s/s'$
\[
D(\omega)=\int_{\langle s'(\omega)\rangle}^{\langle s(\omega)\rangle}
\]
is locally the section at every $\omega$.  
This is well-defined and independent of choices $\LL'$, $s$, and $s'$ by the classical theory of Jacobian varieties, thus we can patch these sections together to get a continuous section on all of $U$.
\end{proof}

\begin{theorem}\label{thm.MtoJ}
$\HH_1(\mathcal M)\cong \HH_1(\SS)^{\tau=-1}$.
\end{theorem}

\begin{proof}
First notice that $\HH_1(\SS)\cong \HH_1(\text{Jac}(\SS))$ naturally, since $\text{Jac}(\SS)$ is defined locally as $\mathbb R^{2\tilde g}$ mod the lattice which is $\HH_1(\SS)$.  This isomorphism is also $\tau$-equivariant, since the $\tau$ action on $\text{Jac}(\SS)$ is defined by lifting from $\SS$.   Thus we need only show that $\HH_1(\mathcal M)\cong \HH_1(\text{Jac}(\SS))^{\tau=-1}$.

By Hitchin (we follow \cite{Hitchin2}, Theorem 8.1 closely), for every $\omega$, there exists some $L_\omega$ such that $L\mapsto L\otimes L_\omega$ is an isomorphism:
\[
\mathcal M|_{\{\omega\}}\to\text{Pic}(S_\omega,M)=\text{Jac}(S_\omega)^{\tau=-1}.
\]
It follows from this that for any $L_\omega'$ such that the degree of $L_\omega'$ equals the degree of $L_\omega$, the map $L\to L\otimes L_\omega'$ induces the same map 
\[
H_1(\mathcal M|_{\{\omega\}},\mathbb Z)\to H_1(\text{Jac}(S_\omega),\mathbb{Z})^{\tau=-1}.
\]
as does $L\to L\otimes L_\omega$, since homology is translation invariant.  Thus there is a well-defined map (as sets):
\[
\HH_1(\mathcal{M})\to \HH_1(\text{Jac}(\SS)).
\]

To any local section of $\mathcal M\to H^s$, we have a map $\Phi:U\times V\to U\times V\otimes K_M$ on vector bundles over $U\times M$.  The vector bundle $U\times V\to U\times M$ has a smooth operator $\bar\partial_\omega$ defining the holomorphic structure on the fibers, $V$.  $\Phi$ has the property that $\det(\Phi(\omega))=\omega$.  $\Phi$ pulls back to $\tilde\Phi$ acting on $\Pi^*(U\times V)$ over $\SS|_U$.  By construction, $\tilde \Phi$ preserves the holomorphic structure on $\Pi^*(U\times V)|_{S_\omega}$.  $\tilde\Phi:\Pi^*(U\times V)\to\Pi^*(U\times V)\otimes \Pi^*(K_M)$.  We also have the canonical section $\sqrt{\text{det}\Phi}$.

Let $\LL_\Phi=\text{ker}(\tilde\Phi-\sqrt{\text{det~}\Phi})$ be the line bundle on $\SS|_U$, which is holomorphic along the $S_\omega$, and $d$ be its degree.  The degree is independent of choice of section of $\mathcal M$.  Choose a line bundle $\LL$ of degree $-d$ by Lemma \ref{lem.linebun}.  Then the map $\LL_\Phi\to \LL_\Phi\otimes \LL$ gives a degree zero line bundle associated to any $\Phi$ which is holomorphic along the $S_\omega$.  By Proposition \ref{prop.transportsec}, we get a map from smooth sections of $\mathcal{M}|_U$ to smooth sections of $\text{Jac}(\SS)|_U$.  This shows that the bundles $\mathcal{M}$ and $\text{Jac}(\SS)$ are locally (non-canonically) isomorphic.  However, any choice of isomorphism gives the same map:
\[
\HH_1(\mathcal{M})|_U\hookrightarrow \HH_1(\text{Jac}(\SS))|_U
\]
defined before, showing that this map is continuous, and thus inducing a global isomorphism to the image, $\HH_1(\text{Jac}(\SS))^{\tau=-1}$.
\end{proof}

\section{Relating to the Burau representation}\label{sec.relating}

In Section \ref{sec.burau}, we defined the Burau representation, while in Theorem \ref{thm.MtoJ}, we showed the homology of the Hitchin bundle $\HH_1(\MM)$ was equivalent to $\HH_1(\SS)^{\tau=-1}$.  In this section, we show that $\HH_1(\SS)^{\tau=-1}$ is indeed a specialization of the Burau bundle.  Let $k=2$.

Recall the space $\WW$ from Section \ref{sec.burau}, letting $n=4g-4$.  This is the space of points in $M^{[4g-4]+1}$ along with a winding number function.  For any $\omega$, we can construct a (mod 2) winding number function $w_\omega$ defined on closed loops $\alpha:[0,1]\to M\setminus \langle \omega\rangle$ by:
\[
w_\omega(\alpha) =\frac{\sqrt{\omega(\alpha(1))}}{\sqrt{\omega(\alpha(0))}} = \pm 1,
\]
which is well-defined, and descends to $H_1(M\setminus  \langle \omega\rangle,\mathbb Z)$.  Also, $w_{c\omega}=w_\omega$ for all nonzero constants, $c$.  Notice that we computed a similar function explicitly in Theorem \ref{thm.separated}.

If $\langle\omega\rangle=X$, then for any $m\notin X$, $\BB|_{(m,X,w_\omega)}$ and $\SS|_\omega$ are isomorphic, yet there are two choices of isomorphism, depending on the image of the basepoint, $\pm\sqrt{\omega(m)}$.  Recall also that the representation 
\[
\pi_1(\WW)\to \text{Aut}(H_1(\BB|_{(mXw)},\mathbb Z))/\langle t\cdot \text{Id}\rangle
\]
is trivial on the fibers of $\WW\to M^{[4g-4]}$.  Therefore the representation descends to a well-defined representation of $\pi_1(M^{[4g-4]})$, and likewise for $\pi_1(PH^s)\to\text{Aut}(H_1(S_{\omega},\mathbb Z))/\langle t\cdot \text{Id}\rangle$. 
This gives us 

\begin{proposition}\label{prop.withdiagram}
Let $\omega\in H^s$.  Under $PH^s\to M^{[4g-4]}$, and either identification, $f:B_{mXw_\omega}\to S_{\omega}$, where we use $\bar\omega$ to denote the element of $PH^s$ to avoid confusion, $f$ induces an intertwining isomorphism, 
\[
\hat f:  \text{Aut}(H_1(\BB|_X,\mathbb Z))/\langle t\cdot \text{Id}\rangle
     \to\text{Aut}(H_1(S_{\omega},\mathbb Z))/\langle t\cdot \text{Id}\rangle
\]
realized by the commutative diagram:
\[
\xymatrix{
\pi_1(PH^s)\ar[r]\ar@{->}[d]& \text{Aut}(H_1(S_{\omega},\mathbb Z))/\langle t\cdot \text{Id}\rangle\\
\pi_1(M^{[4g-4]})\ar[r]& \text{Aut}(H_1(\BB|_X,\mathbb Z))/\langle t\cdot \text{Id}\rangle\ar[u]_{\hat f}
}
\]
\end{proposition}

We will prove this proposition in a moment.  However, first realize an important corollary:

\begin{corollary}\label{cor.1tau}
In the representation, $\pi_1(H^s)\to\text{Aut}(H_1(S|_{\omega},\mathbb Z))$, the kernel 
\[
\text{ker}(\pi_1(H^s)\to \pi_1(M^{[4g-4]}))
\]
acts by $\{1,\tau\}$.
\end{corollary}


\begin{proof}[Proof of Proposition \ref{prop.withdiagram}]
First note that $(\hat f\sigma)(\alpha)=(f^{-1})^*(\sigma(f^*\alpha))$ is an isomorphism, as it is invertible.  Also, if $f'=tf$, then $\hat{f'}(\gamma)=(t\hat f t)\gamma$, which gives a different diagram, however, one commutes if and only if the other commutes, since $t\hat f t$ agrees with $\hat f$ on the image of $\pi_1(H^s)$, since $t$ commutes with the monodromy action as observed at the end of Section \ref{sec.burau}.

Let $U\subset M$ be contractible such that the canonical bundle is trivializable on $U$ with $\mathbb C^*$-equivariant trivialization $\bar\phi:K_M|_U\to U\times \mathbb C$.  This induces a trivialization $\phi:K_M^2|_U\to U\times\mathbb C$, respecting tensor product ($\phi\cong\bar\phi\otimes_U\bar\phi$).  Consider the space 
\[
A=\{(z,m,X)\in\mathbb C^*\times U\times PH^s|m\notin X\}.
\]
As $U$ has infinite cardinality (thus greater than $4g-4$), $A$ has a surjective map $A\twoheadrightarrow PH^s$.  In fact, for every $(z,m,X)\in A$, there exists a unique $\omega\in H^s$ such that $\omega|_X=0$ and $\phi(\omega(m))=(m,z^2)$.  Thus we get a lift to $H^s$ which is surjective as well (by $\mathbb C^*$ equivariance, say).  

This gives rise to a commutative diagram:

\[
\xymatrix{
\HH_1(\mathcal S\ar[d]) & a^*\HH_1(\mathcal S)\ar[l]\ar[d]  &  b^*\HH_1(\BB)\ar[r]\ar[d] & \HH_1(\BB)\ar[d] \\
H^s\ar[d] & A\ar@{->>}[l]_{a}\ar@{->>}[dl] \ar@{=}[r] & A\ar[r]^{b}\ar@{->>}[dr] & \WW \ar[d]\\
PH^s\ar@{^(->}[rrr] & & & M^{[4g-4]} \\
}
\]

We intend to show that $a^*\HH_1(\mathcal S)$ and $b^*\HH_1(\BB)$ are isomorphic.  Once we have this, clearly 
\[
a^*\HH_1(\mathcal S)/\langle t\cdot \text{Id}\rangle\cong b^*\HH_1(\mathcal B)/\langle t\cdot \text{Id}\rangle
\]
are trivial on fibers of $A\to M^{[4g-4]}$, so the result will follow.

We in fact will show something stronger, that $a^*\SS$ and $b^*\BB$ are isomorphic.  This isomorphism boils down to a universality statement, Lemma \ref{lem.universal}, that there is only one double cover with given data $(X,w,m)$ of ramification divisor, winding number and marked point in the cover.

Let $\WW'=bA$.  Then $b^*\HH_1(\BB)\cong \HH_1(\BB)|_{\WW'}\times \mathbb C^*$.  Let 
\[
Q=\{
((z,m,X),b)\in\HH_1(\SS)|z=1
\}.
\]
Then $a^*\HH_1(\SS)\cong Q\times \mathbb C^*$ by the map 
\[
((z,m,X),b)\mapsto ((1,m,X),\frac bz,z).
\]
This map is well-defined since if $\omega(\pi_\omega(b))=b^2$, then $\frac\omega{z^2}(\pi_{\frac\omega{z^2}}(\frac bz))=\frac{b^2}{z^2}\in S_{\frac\omega{z^2}}$.  It is also clearly continuous with continuous inverse.  Now applying the Lemma \ref{lem.universal} to $Q\to \WW'$ and $\HH_1(\BB)|_{\WW'}$, we get $Q\cong \HH_1(\BB)|_{\WW'}$, so
\[
a^*\HH_1(\SS)\cong Q\times \mathbb C^* \cong \HH_1(\BB)|_{\WW'}\times \mathbb C^*\cong b^*\HH_1(\BB).
\]

\end{proof}

\section{Combinatorics of $\HH_1(\MM)$}

In this section, we use the results of Section \ref{sec.prym} to turn $\HH_1(\MM)$ into combinatorial data.  First notice
\begin{align}
H_1(\MM|_{\{\omega\}})&\cong H_1(\text{Prym}(S_\omega,M),\mathbb Z)\notag\\
&\cong\mathbb RE/(\mathbb RF+\frac12\mathbb ZE).\notag
\end{align}
Also recall that the monodromy of $\HH_1(\MM)\to H^s$ is modelled on $\HH_1(\SS)$, so that studying the effect on $\pi_1(H^s)$ on $\frac12\mathbb ZE$ gives us this by 
\[
\HH_1(\SS)|_{\{\omega\}}\cong\frac12\mathbb ZE/(\mathbb RF\cap\frac12\mathbb ZE).
\]
Given Theorem \ref{thm.iso}, then, our only remaining task is to compute for edges $e_\alpha, e_\beta\in E$, the values
\[
\sigma_{e_\alpha}\frac12\phi(e_\beta)\in\frac12\mathbb Z\phi(E).
\]

First, however, we discuss a result making the introduction of the Burau representation more salient.  Let $B_{\Lambda_2}$ be the Burau module.  This is a free module over $\Lambda_2=\mathbb Z[t]/\langle t^2-1\rangle$, so 
\[
B_{\mathbb Z}^-=B_{\Lambda_2}/\langle t=-1\rangle
\]
is an honest $\mathbb Z$-module with a ($\mathbb Z$-projective) representation
\[
K=\ker(\pi_1(M^{[n]})\to H_1(M))\to \text{Aut}(B_{\mathbb Z}^-)/\langle \pm1\rangle.
\]
We choose this notation based on the observation that 
\[
B_{\mathbb Z}=B_{\Lambda_1}=B_{\Lambda_2}/\langle t=+1\rangle.
\]

Notice that in $H_1(\text{Prym}(S_\omega,M),\mathbb Z)$, $\tau$ acts by $-1$, so 
\[
\pi_1(H^s)\to\text{Aut}(H_1(\text{Prym}(S_\omega,M),\mathbb Z)
\]
descends to a projective representation
\[
\pi_1(PH^s)\to\text{Aut}(H_1(\text{Prym}(S_\omega,M),\mathbb Z)/\langle \pm1\rangle.
\]

\begin{theorem}
The standard homomorphism $H_1(\text{Prym}(S_\omega,M),\mathbb Z)\to H^1(S_\omega,\mathbb Z)\to B_{\mathbb Z}^-$ induces a commutative diagram and $\pi_1(PH^s)$-module isomorphism:
\[
\xymatrix{
\pi_1(PH^s)\ar[r]\ar@{->}[d]& \text{Aut}(H_1(\text{Prym}(S_\omega,M),\mathbb Z))/\langle \pm1\rangle\\
K\ar[r]& \text{Aut}(B_{\mathbb Z}^-)/\langle \pm 1\rangle\ar[u]_{\cong}
}
\]
\end{theorem}

\begin{proof}
In fact, commutativity of this diagram follows immediately from Proposition \ref{prop.withdiagram}.  The crux of the statement is that these modules are isomorphic.  However, this also follows from Proposition \ref{prop.withdiagram}, once one observes that $\hat f$ is an intertwining isomorphism and is $\Lambda_2$-equivariant.
\end{proof}

\section{monodromy}

Now we are ready to study the monodromy action.  We know this action is generated by the transpositions, so we let $\bar\sigma\in\rho_*\pi_1(H^s)$ realize one of these transpositions.  By Corollary \ref{cor.1tau} we know that if we choose any lift $\sigma\in\pi_1(H^s)$, then any other lift acts in the manner of $\sigma$ or $\tau\sigma$.  The degree of freedom comes from the ambiguity of identifying $S_\omega$ with $S_{\sigma\omega}$, as we have discussed.  

In this section, we study the effect of transpositions on $M$.  Recall that $\Gamma$ makes $M$ a polyhedron, $c:S\to M$ is a double cover ramified at marked points, one per face of $M$.  We would like to study the action of $\rho_*\pi_1(H^s)$ on the edge set of $\cbG$.  The faces of $M$ are labelled by vertices of $\cbG$.  Let $X_j$ be a face of $\Gamma$ and let $\{x_j\}= X_j\cap X$ be the corresponding vertex of $\check\Gamma$.  Assume $x_1$ and $x_2$ are neighbors across $e_\alpha$, with associated transposition $\sigma_\alpha$.  Assume, once and for all, that every transposition is realized by each point $x_j$ moving counterclockwise around $e$.   

\begin{theorem}\label{thm.repwithint}
If $e_\alpha\in E$ and $\overline\sigma_\alpha\in\rho_*\pi_1(H^s)$, $\overline\sigma_\alpha$ has a lift $\sigma_\alpha\in\pi_1(H^s)$ such that for any $e\in\mathbb{Z}E$, 
\[
\sigma_\alpha e=e-(e\cdot e_\alpha)e_\alpha.
\]
\end{theorem}

\begin{proof}
We choose the lift $\sigma_\alpha$  of $\bar\sigma_\alpha$ which acts homotopically trivially on $S\setminus c^{-1}(X_1\cup X_2)$.  Edges away from $X_1\cup X_2$ are unaffected by this transposition, however any edge from $x_1$ or $x_2$ may be altered.  $X_1\cup X_2$ is simply connected by construction, and has two ramification points.  Recall that we are studying the effect of $\sigma$ on $E\subset \bar E$.  Notice that the loop $\phi e_\alpha$ is in fact homotopic to one component of $\partial c^{-1}(X_1\cup X_2)$, thus in fact to a loop lying outside $c^{-1}(X_1\cup X_2)$.  This implies that $\sigma_\alpha$ acts trivially on $e_\alpha$.  In fact, we have seen that $\sigma_\alpha$ acts trivially on any $e\in E$ such that $e\cdot e_\alpha=0$, where $\cdot$ is the intersection pairing on $S$.

Notice that if $e_1\cdot e_\alpha=e_2\cdot e_\alpha=+1$, then $e_1-e_2$ is homotopy equivalent to a curve supported away from $c^{-1}(X_1\cap X_2)$, so $\sigma e_1-e_1=\sigma e_2-e_2$.  By linearity, then,
\[
\sigma_\alpha e =e-(e\cdot e_\alpha)e_q
\]
for some as yet undetermined $e_q$.

The set $c^{-1}(X_1\cup X_2)$ may be realized topologically as the union of two annuli joined along their inner circles (preserving orientation).  The effect of $\sigma$ on each annulus is a positive half-twist of these circles.  Whe rejoined, we see that $\sigma$ is the operation of a Dehn twist at $e_\alpha$.  Thus 
\[
\sigma_\alpha e =e\pm(e\cdot e_\alpha)e_\alpha.
\]

Notice that 
\[
e\cdot (\sigma_\alpha e-e) =\pm(e\cdot e_\alpha)^2,
\]
and it is enough to show that there is some $e$ such that that this quantity is negative.

In some small neighborhood of $e_\alpha$, the monodromy is computed as in the Figures 24 and 25.  $e$ is a curve which crosses the branch cut $e_\alpha$, so the dark curves lie on the upper branch and the light curves lie on the lower branch.  Applying the monodromy, we get Figure 24.

\begin{figure}
\noindent
\includegraphics[width=1\textwidth]{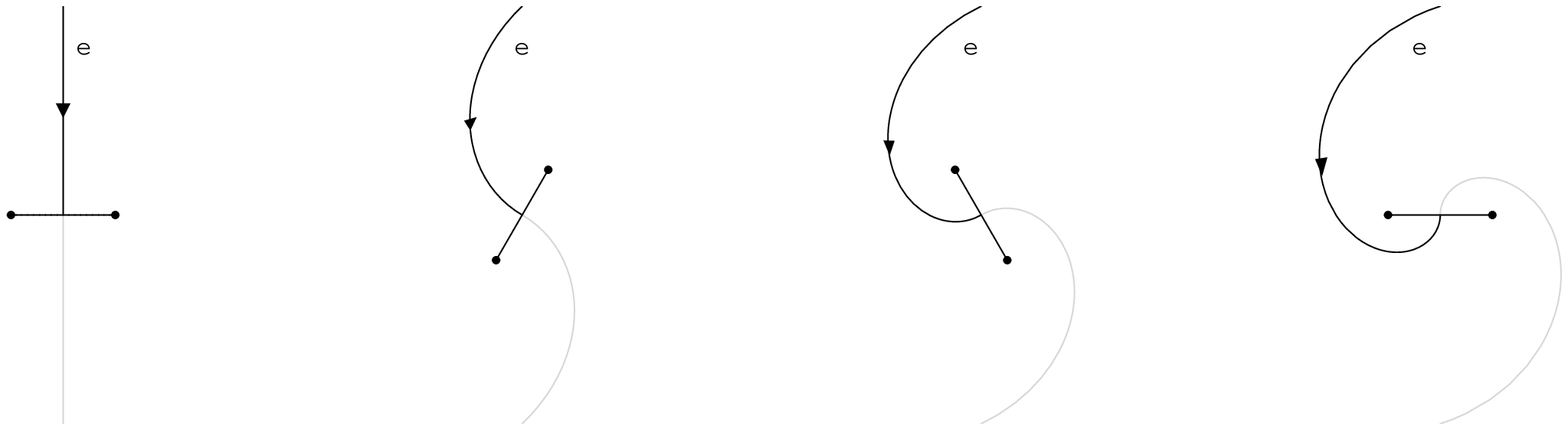}

Figure 24. Monodromy action of $\sigma_\alpha$ on $e$.
\end{figure}

One then graphs $e$ and $ \sigma_\alpha e-e$ as in Figures 25 and 26, and finds that their intersection is $-1$ as desired.  Of course the sign was originally chosen by the choice of transposition $\sigma_\alpha^{\pm}$.

\begin{figure}
\noindent
\includegraphics[width=.3\textwidth]{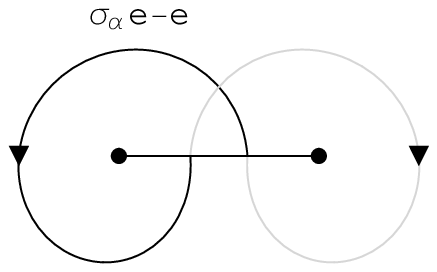}

Figure 25. Computing $\sigma_\alpha e-e$.
\end{figure}

\begin{figure}
\noindent
\includegraphics[width=.3\textwidth]{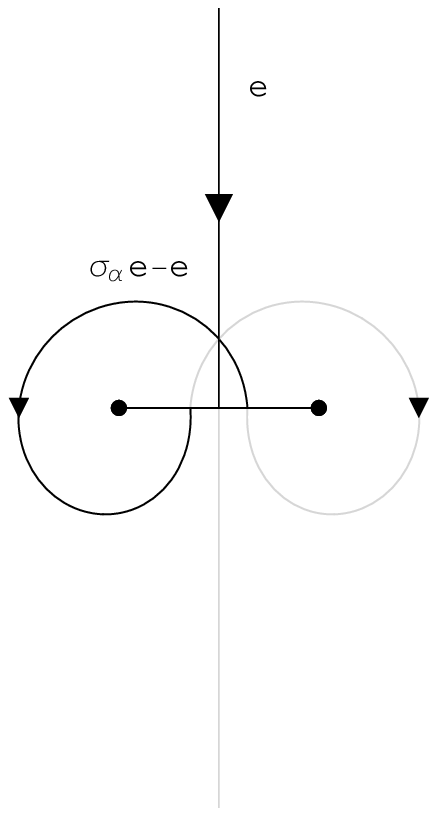}

Figure 26. Computing $e\cdot(\sigma_\alpha e-e)$.
\end{figure}

\end{proof}

\section{A note regarding non-hyperelliptic curves}\label{sec.generalize}

In this section, we show how, given Theorem \ref{thm.iso}, for general $M$, that the rest of the work in this paper still applies.  Indeed, we need only that there is a compatible cellular decomposition:

\begin{theorem}\label{thm.generalM}
Let $N$ be any curve of genus $g\geq 3$ with $\eta\in H^0(N,K_N^2)^s$, and $M$ hyperelliptic with $\omega\in H^0(M,K_M^2)^s$ as in Section \ref{sec.g3}.  Then $S_\eta$ and $S_\omega$ are homeomorphic as double covers.  
\end{theorem}

\begin{proof}
Let $X=\langle\eta\rangle$ and choose $x_0\in N\setminus X$.  Any closed curve through $x_0$ missing $X$ is homotopic in $N$ to another such curve with winding number zero (by dragging the original curve past some subset of $X$).  Thus there is a set of closed curves $a_j$, $b_j$ of winding number $0$ such that $[a_j], [b_j]$ is a symplectic basis for $H_1(N,\mathbb Z)$.  Excising these curves, we find a contractible open set $U_N\subset N$ containing all of the ramification points, $X$.  By construction, this set lifts to a two-component set in $S_\eta$.  Also, $U$ is the interior of a polygon, from which one can form $N$ by identifying sides in a standard way, and $S_\eta$ is defined uniquely by this polygon.

We used no data about $N$ in constructing this decomposition, so we could do the same for $M$ and $\omega$, finding $U_M$.  However, $U_M$ and $U_N$ are homeomorphic as punctured open disks, and indeed, there is a homeomorphism preserving the gluing relation.  Thus we get two choices for lifts to homeomorphisms of $S_\eta$ and $S_\omega$.
\end{proof}

Now that these are homeomorphic, we may transport the cellular decomposition of $(S_\omega,M)$ to $(S_\eta,N)$.  This is the only other data that relies on hyperellipticity, so we see that, indeed, if Theorem \ref{thm.iso} holds in general, then all other arguments hold as well.

\section{acknowledgements}
I'd like to thank my advisor, V. Ginzburg, for helpful suggestions and discussions.

\small{

}

\end{document}